\documentclass[a4paper,11pt]{amsart}
\usepackage{amsmath}
\usepackage{amsthm}
\usepackage{amssymb}
\usepackage{amscd}
\usepackage{amsfonts} 
\usepackage[all]{xy}

\usepackage{ascmac}

\title{The spinor and tensor fields with higher spin on spaces of constant curvature}
\author{Yasushi Homma }
\address{Department of Mathematics, Faculty of science and engineering, Waseda University, 3-4-1 Ohkubo, Shinjuku-ku, Tokyo, 169-8555, JAPAN.}
  \email{homma\_yasushi@waseda.jp} 
\author{Takuma Tomihisa }
\address{Department of Pure and applied Mathematics, Graduate school of fundamental science and engineering Waseda University, 
  3-4-1 Ohkubo, Shinjuku-ku, Tokyo, 169-8555, JAPAN.}
  \email{taku-tomihisa@akane.waseda.jp} 
  
\keywords{the Lichnerowicz Laplacian, the (higher spin) Dirac operator, Weitzenb\"ock formulas, generalized gradients, space of constant curvature, harmonic analysis on spheres }
  \subjclass[2010]{53C27, 58C40, 58J50}
   \date{}

\theoremstyle{plain}
\newtheorem{theorem}{Theorem}
\newtheorem{proposition}[theorem]{Proposition}
\newtheorem{corollary}[theorem]{Corollary}
\newtheorem{lemma}[theorem]{Lemma}
\theoremstyle{definition}

\numberwithin{equation}{section}

\theoremstyle{remark}
\newtheorem{rem}{Remark}[section]

\newtheorem{exam}{Example}[section]

\newcommand{\Hom}{\mathrm{Hom}}

\newcommand{\id}{\mathrm{id}}

\newcommand{\la}{\langle}
\newcommand{\ra}{\rangle}

\newcommand{\SO}{\mathrm{SO}}
\newcommand{\Spin}{\mathrm{Spin}}
\newcommand{\Scal}{\mathrm{Scal}}
\newcommand{\Sym}{\mathrm{Sym}}
\allowdisplaybreaks[1]
\begin{document}

\maketitle

\begin{abstract}
In this article, we give all the Weitzenb\"ock-type formulas among the geometric first order differential operators on the spinor fields with spin $j+1/2$ over Riemannian spin manifolds of constant curvature. Then we find an explicit factorization formula of the Laplace operator raised to the power $j+1$ and understand how the spinor fields with spin $j+1/2$ are related to the spinors with lower spin. As an application, we calculate the spectra of the operators on the standard sphere and clarify the relation among the spinors from the viewpoint of representation theory. Next we study the case of trace-free symmetric tensor fields with an application to Killing tensor fields. Lastly we discuss the spinor fields coupled with differential forms and give a kind of Hodge-de Rham decomposition on spaces of constant curvature. 
\end{abstract}

\tableofcontents

\section{Introduction}\label{sec:1}
Recently many researchers in mathematics and physics have tried to understand geometry and analysis for the spinor fields with higher spin and the symmetric tensor fields. In fact, there are a lot of articles for the Rarita-Schwinger fields with spin $3/2$ and Killing tensor fields, which are special fields satisfying geometric differential equations. One of the methods to define such meaningful geometric differential equations is to use the generalized gradients or the Stein-Weiss operators. The operators are first-order differential operators naturally defined on spinor and tensor fields. In \cite{Fegan}, Fegan classified them and showed they are conformally covariant. In \cite{Br2}, Branson studied the ellipticity and Weitzenb\"ock formulas for generalized gradients using spectra of them on the standard sphere. He also studied another type of Weitzenb\"ock formula in \cite{Br3}. After a decade, the first author of this paper showed an explicit method to construct the Weitzenb\"ock formulas for the generalized gradients in \cite{H1}, which produce a lot of applications, vanishing theorem, eigenvalue estimates and so on. There are also many articles to study the generalized gradients and their applications to mathematics and physics (\cite{Her}, \cite{H2},  \cite{HS}, \cite{Pilca}, etc.). Moving on to analysis, we know that one of the main topics in Clifford analysis is to generalize spherical harmonic analysis on Euclidean space to such spinor and tensor fields. Polynomial solutions called monogenic function, fundamental solutions and Fueter theorem for the higher spin Dirac operators have been constructed in \cite{BSSL1}, \cite{BSSL2}, \cite{ES}, etc. The key to give the fundamental solution is a so-called factorization formula. The $j+1$-st power of the Laplace operator $\Delta$ on $\mathbb{R}^n$ is factorized as $\Delta^{j+1}=D_j\circ A_{2j+1}$, where $D_j$ is a generalized gradient (the higher spin Dirac operator) and $A_{2j+1}$ is a differential operator with order $2j+1$. To consider such formula on curved manifold would be important to analyze higher spin fields and give a variety of applications to geometry. 

 In this article, we combine Clifford analysis and differential geometry. We study the generalized gradients on spinor fields with spin $j+1/2$ over a Riemannian spin manifold $(M,g)$, and try to give a factorization formula for the exponentiation  of the Laplace operator $\Delta$. Since Weitzenb\"ock formulas are too complicated for $j\ge 2$ to handle them on curved manifolds, we assume that $(M,g)$ has a constant sectional curvature. Then we see any operators commute with each other in a way and construct a factorization formula explicitly in Theorem \ref{thm:fact}. To clarify the meaning of the factorization, we show how the spinor fields with spin $j+1/2$ are influenced by the spinor fields with lower spin in Theorem \ref{thm:7}. Remark that, for the case of $j=1$, we need only the assumption of Einstein manifold and can develop fruitful geometry and analysis in \cite{BM}, \cite{HS} and \cite{HT}. Next, we study harmonic analysis on spinor fields with spin $j+1/2$ on the standard sphere as a model case for spinor analysis on a curved space in Section \ref{sec:3}. Interestingly the factorization formula and Weitzenb\"ock formulas yield all the eigenvalues of the generalized gradients from the eigenvalues of the Laplacian, which are calculated by Freudenthal's formula for Casimir operator (cf. \cite{Br1}). This method would be easier to understand than in \cite{Br4} or originally in \cite{BOO}. We also show how the spaces of the spinor fields as $\Spin$-modules relate to each other through the operators. In Section \ref{sec:4}, we discuss the trace-free symmetric tensor fields on a space of constant curvature and give a factorization formula in Theorem \ref{thm:fact2}. We study harmonic analysis on such fields over the standard sphere in Section \ref{sec:5}. As an application, we give a decomposition of the Killing tensor fields on the sphere from the viewpoint of representation theory. In Section \ref{sec:6} and \ref{sec:7}, we discuss the spinor fields coupled with differential forms on a space of constant curvature and harmonic analysis on the sphere. In particular, we give a kind of Hodge-de Rham decomposition for spinor fields with differential forms. In Appendix \ref{app:A}, we show how to calculate the Weitzenb\"ock formulas needed in Section \ref{sec:2}.

\noindent
{\em Acknowledgments.} This work was partially supported by JSPS KAKENHI Grant Number JP19K03480.

\section{Higher spin Dirac operators on spinor fields with spin $j+1/2$}\label{sec:2}

Let $(M,g)$ be an $n$-dimensional Riemannian spin manifold with a spin structure $\Spin(M)$, which is a principal $\Spin(n)$-bundle over $M$ and doubly covering the orthonormal frame bundle $\SO(M)$. Throughout this paper, we assume $n\ge 3$. The case of dimension two is left as an exercise to the reader. We consider the spin $j+1/2$ (unitary) representation $\pi_j$ on $W_j$ of $\Spin(n)$ for $j=0,1,\cdots$. For example, $(\pi_0,W_0)$ is a usual spinor representation. For $n=2m$, the space $W_j$ is decomposed into the direct sum of $W_j^+$ and $W_j^-$. Each $W_j^{\pm}$ is an irreducible $\Spin(n)$-module whose highest weight is
\[
(j+1/2,\underbrace{1/2,\cdots,1/2}_{m-2},\pm 1/2)=(j+1/2,(1/2)_{m-2},\pm 1/2).
\]
Here $(1/2)_k$ denotes a sequence $1/2,\cdots,1/2$ with length $k$ as an abbreviation. We write the direct sum representation $\pi_j^+\oplus \pi_j^-$ by $\pi_j$. For $n=2m-1$, the spin $j+1/2$ spinor space $W_j$ is an irreducible $\Spin(n)$-module whose highest weight is $(j+1/2,(1/2)_{m-2})$. It follows from Weyl's dimension formula that the dimension of $W_j$ is
\begin{gather*}
\dim W_j=2^{[n/2]}\binom{n+j-2}{j}, \quad \dim W_j^{\pm}=2^{[n/2]-1}\binom{n+j-2}{j}.
\end{gather*}
The representation $\pi_j$ (resp. $\pi_j^{\pm}$) induces a vector bundle $S_j$  (resp. $S_j^{\pm}$) associated to the principal bundle $\Spin(M)$, 
\[
S_j:=\Spin(M)\times_{\Spin(n)}W_j.
\]
We call a section of $S_j$ {\it a spin $j+1/2$ field} or {\it a spinor field with spin $j+1/2$}. 

From now on, we study some basic properties for the first order differential operators naturally defined on the space of the spin $j+1/2$ fields $\Gamma(S_j)$. They are called {\it generalized gradients} or {\it Stein-Weiss operators}, which are defined by composing the orthogonal bundle projections and the covariant derivative $\nabla$ induced by the Levi-Civita connection, \cite{Br2}, \cite{Fegan}. The covariant derivative on $\Gamma(S_j)$ is 
\[
\nabla:\Gamma(S_j)\ni \phi \mapsto \nabla \phi=\sum \nabla_{e_i}\phi \otimes e_i\in \Gamma(S_j\otimes TM^c),
\]
where $\{e_i\}_i$ is a local orthonormal frame and $TM^c$ is $TM\otimes \mathbb{C}\cong T^{\ast}M\otimes \mathbb{C}$ by Riemannian metric $g$. We split the fiber $W_j\otimes \mathbb{C}^n$ into the sum of $\Spin(n)$-modules, 
\[
W_j\otimes \mathbb{C}^n=W_{j+1}\oplus W_{j,1}\oplus W_j\oplus W_{j-1},
\]
where $W_{j,1}$ is a $\Spin(n)$-module with the highest weight
\[
\begin{cases}
(j+1/2,3/2,(1/2)_{m-2})\oplus (j+1/2,3/2,(1/2)_{m-3},-1/2) & \textrm{for $n=2m$}, \\
(j+1/2,3/2,(1/2)_{m-3})& \textrm{for $n=2m-1$}.
\end{cases}
\]
We note that $W_{j,1}$ doesn't appear for $n=3$ or $j=0$, and $W_{j-1}$ doesn't appear for $j=0$. 

The space $W_j$ has a $\Spin(n)$-invariant Hermitian inner product (unique up to a constant factor), so that the above decomposition is orthogonal. Since the fiber metric on $S_j$ is induced from the inner product, we have the orthogonal bundle projection onto each irreducible summand. For example, composing $\nabla$ and the projection 
\[
\Pi_j:S_j\otimes TM^c\to S_j,
\]
we have so-called {\it the higher spin Dirac operator},
\begin{equation}
\tilde{D}_j:=\Pi_j\circ \nabla,\quad \Gamma(S_j)\xrightarrow{\nabla}\Gamma(S_j\otimes TM^{c})\xrightarrow{\Pi_j}\Gamma(S_j). \label{h-Dirac}
\end{equation}
In this manner we construct four generalized gradients on $\Gamma(S_j)$ and name them as follows;
\begin{align*}
& \tilde{T}_j^+:\Gamma(S_j)\to \Gamma(S_{j+1}) & \textit{the (first) twistor operator},\\
 & U_j:\Gamma(S_j)\to \Gamma(S_{j,1}) & \textit{the (second) twistor operator},\\
 & \tilde{D}_j:\Gamma(S_j)\to \Gamma(S_j) & \textit{the higher spin Dirac opeator},\\
 & \tilde{T}_j^-:\Gamma(S_j)\to \Gamma(S_{j-1}) & \textit{the co-twistor opeator}.
\end{align*}
Here we set $U_j=0$ for $n=3$, $U_0=0$ and $\tilde{T}_0^-=0$. When $n=2m$, each operator has the form of $2\times 2$ matrix as
\[
\tilde{D}_j=\begin{pmatrix}
 0 & \ast\\
\ast & 0 
\end{pmatrix},\quad \tilde{T}_j^+,U_j,\tilde{T}_j^-=\begin{pmatrix}
\ast &  0\\
 0 &  \ast
\end{pmatrix}
\]
along the decomposition $\Gamma(S_j)=\Gamma(S_j^+)\oplus \Gamma(S_j^-)$. We introduce a $L^2$-inner product on $\Gamma(S_j)$ by
\[
(\phi_1,\phi_2):=\int_M \la\phi_1,\phi_2\ra \mathrm{vol}_g,\quad (\phi_1,\phi_2\in \Gamma(S_j)).
\]
Then the co-twistor operator $\tilde{T}_j^-$ is a non-zero constant multiple of the formal adjoint of the twistor operator $\tilde{T}_{j-1}^+$ with spin $j-1/2$,
\[
(\tilde{T}_{j-1}^+)^{\ast}=c\tilde{T}_j^-:\Gamma(S_{j})\to \Gamma(S_{j-1}). 
\]

In the case of $j=0$, the operator $\tilde{D}_0$ seems to coincide with the Dirac operator. However we know $\tilde{D}_0=1/\sqrt{n}D$. Therefore we need a normalization of $\tilde{D}_j$ to get the (usual) higher spin Dirac operator $D_j$ in other articles (\cite{BSSL2} etc.). An advantage of our definition with projection is that we can apply Weitzenb\"ock formulas in \cite{H1} and \cite{H2}. Furthermore such formulas allow us to give a normalization of $\tilde{D}_j$ explicitly. Before normalizing the operators, we show some known results for analytic properties of the generalized gradients. 
\begin{proposition}\label{prop:1}
\begin{enumerate}
\item $\tilde{D}_j,\tilde{T}_j^{\pm},U_j$ are conformally covariant. 
\item $\tilde{D}_j$ is a first order, elliptic and (formally) self-adjoint differential operator. 
\item $\tilde{T}_{j}^+$ is an overdetermined elliptic operator in the sense that the principal symbol $\sigma_{\xi,x}(\tilde{T}_j^+)$ is injective for every $x$ in $M$ and non-zero covector $\xi$ in $T^{\ast}_xM$. Then there is an orthogonal decomposition for the space of sections on a compact Riemannian spin manifold $M$,
\[
\Gamma(S_j)=\ker (\tilde{T}_{j-1}^+)^{\ast} \oplus \tilde{T}_{j-1}^+(\Gamma(S_{j-1}))=\ker \tilde{T}_{j}^- \oplus \tilde{T}_{j-1}^+(\Gamma(S_{j-1})).
\]
Moreover $\tilde{T}_j^+$ is a differential operator of finite type in the sense that the kernel of $\tilde{T}_j^+$ is finite dimensional (even if $M$ is non-compact). 
\end{enumerate}
\end{proposition}
\begin{proof}
The first claim follows from the fact that generalized gradient is conformally covariant \cite{Fegan}, \cite{H1}. As for the second claim, a general result for the ellipticity is known in \cite{Br2}, \cite{Pilca}. We can also get the ellipticity for $D_j$ from Remark \ref{symbol}. For the third one, we require a discussion. Let $U$ and $V$ be irreducible $\Spin(n)$-modules and $W$ be the Cartan summand (the top irreducible summand) of $U\otimes V$. We consider the orthogonal projection $\Pi$ onto an irreducible summand of $U\otimes V$. Then it is shown in \cite{KPOWZ} that the following two propositions are equivalent: 
\begin{itemize}
\item for $u\otimes v$ in $U\otimes V$, $u\otimes v\neq 0$ $\Rightarrow$ $\Pi(u\otimes v)\neq 0$. 
\item $\Pi$ is the projection onto $W$.
\end{itemize}
Applying this to $U=W_j$, $V=\mathbb{C}^n=T_xM^c$ and $W=W_{j+1}$, we know $\Pi_{j+1}(u\otimes \xi)\neq 0$  for non-zero $u\otimes \xi$ in $W_j\otimes \mathbb{C}^n$.  On the other hand, there is a result in \cite{SW2004} for finite type of differential equations such that this condition holds if and only if the generalized gradient is finite type. Thus, we have proved $\tilde{T}_j^+$ is finite type and $\sigma_{\xi,x}(\tilde{T}_j^+)$ is injective. The decomposition of the space of sections with respect to a overdetermined elliptic operator is well-known in \cite{Be}. 
\end{proof}
\begin{rem}
Let $S_{\rho}$ be an irreducible vector bundle with the highest weight $\rho=(\rho^1,\cdots,\rho^m)$. Then we have the generalized gradient $D_{\rho+\mathbf{e}_1}^{\rho}$ from $S_{\rho}$ to the Cartan summand $S_{\rho+\mathbf{e}_1}$ in $S_{\rho}\otimes TM^c$ where $\rho+\mathbf{e}_1=(\rho^1+1,\rho^2,\dots, \rho^m)$. In the same manner, the operator $D_{\rho+\mathbf{e}_1}^{\rho}$ is only the generalized gradient of finite type on $S_{\rho}$.
\end{rem}
We shall give another definition of the higher spin Dirac operator by ``twisting'' the spinor bundle $S_0$ with symmetric tensor bundle in \cite{BSSL2}. Let $\mathrm{Sym}^j=\mathrm{Sym}^j(TM^c)$ be the $j$-th symmetric tensor product bundle of $TM^c$ over $(M,g)$ and $\mathrm{Sym}^j_0$ be its primitive irreducible component whose fiber is an irreducible $\Spin(n)$-module with the highest weight $(j,0_{m-1})$. The bundle $\mathrm{Sym}^j$ is the direct sum of $\mathrm{Sym}_0^j$ and lower parts, 
\begin{equation}
\mathrm{Sym}^j= \mathrm{Sym}_0^j\oplus g\cdot \mathrm{Sym}_0^{j-2}\oplus g^2\cdot \mathrm{Sym}_0^{j-4}\oplus \cdots \oplus g^{[n/2]}\cdot \begin{cases} TM^c & \textrm{$j$ odd} \\
 M\times \mathbb{C}& \textrm{$j$ even},
\end{cases}\label{symj}
\end{equation}
where $g^k\cdot$ is the symmetric tensor product of $g^k=g\cdots g$. We consider the tensor bundle $S_0\otimes \mathrm{Sym}^j_0$ and the twisted Dirac operator defined by
\[
D(j)=\sum_{k=1}^n (e_k\cdot \otimes \id_{\mathrm{Sym}^j_0})\circ \nabla_{e_k}.
\]
Here $\nabla$ is the covariant derivative on $S_0\otimes \mathrm{Sym}^j_0$, and $e_k\cdot$ is the Clifford multiplication by $e_k$. Along the decomposition of the bundle $S_0\otimes \mathrm{Sym}^j_0=S_j\oplus S_{j-1}$, we can realize $D(j)$ as a $2\times 2$ matrix. As stated in \cite{H2}, each component is a generalized gradient up to a nonzero multiplicative constant, 
\begin{equation}
D(j)=\begin{pmatrix}
D_j & T_{j-1}^+ \\
T_j^- & D_{j-1}'
\end{pmatrix},\qquad \Gamma(S_j\oplus S_{j-1})\to \Gamma(S_j\oplus S_{j-1}).
\label{twisted-Dirac}
\end{equation}
Since $D(j)$ is a formally self adjoint operator, we get 
\[
D_j^{\ast}=D_j,\quad (D_{j-1}')^{\ast}=D_{j-1}',\quad (T_{j-1}^+)^{\ast}=T_j^-.
\]
We shall see the square of $D(j)$ gives Weitzenb\"ock formulas. From twisted Lichnerowicz formula,
\begin{equation}
D(j)^2=\Delta_{S_0\otimes \mathrm{Sym}^j_0}+\frac{\Scal}{8}-\frac{1}{2}\id_{S_0}\otimes R_{\mathrm{Sym}^j_0}.
\label{eqn:square1}
\end{equation}
First, we explain the third term on the right side, the curvature action $R_{\mathrm{Sym}^j_0}$ on $\mathrm{Sym}^j_0(TM^c)$. Let $\{e_i=[p,\mathbf{e}_i]\}_{i=1}^n$ be a local orthonormal tangent frame of $TM=\SO(M)\times_{\SO(n)}\mathbb{R}^n$ and $\{e_{ij}=e_i\wedge e_j\}_{1\le i<j\le n}$ be a local frame of 
\[
\Lambda^2(T^{\ast}M)\simeq \mathfrak{o}(TM)=\mathfrak{so}(TM)=\SO(M)\times_{\SO(n)}\mathfrak{so}(n).
\]
We have an action of local section $e_{ij}$ of $\mathfrak{so}(TM)$ on an associated vector bundle $S_{\rho}:=\Spin(M)\times_{\Spin(n)}W_{\rho}$ with respect to a (not necessarily irreducible) representation $(\pi_{\rho},W_{\rho})$:
\[
\mathfrak{so}(TM)\times S_{\rho}\ni (e_{ij}=[p,\mathbf{e}_i\wedge \mathbf{e}_j],\phi=[p,v])\mapsto \pi_{\rho}(e_{ij})\phi:=[p,\pi_{\rho}(\mathbf{e}_i\wedge \mathbf{e}_j)v]\in S_{\rho}.
\]
Then we define the curvature action $R_{\rho}$ on $S_{\rho}$ by
\[
\begin{split}
R_{\rho}&= \frac{1}{2}\sum_{1\le i,j,k,l\le n} R_{ijkl}\pi_{\rho}(e_{ij}e_{kl})\\
&=\frac{1}{2}\sum_{1\le i,j,k,l\le n} W_{ijkl}\pi_{\rho}(e_{ij}e_{kl})
+2\sum_{1\le i,k,j\le n} E_{ij}\pi_{\rho}(e_{ik}e_{jk})+\frac{\pi_{\rho}(c_2)}{n(n-1)}\Scal\\
&=:\mathrm{Weyl}_{\rho}+\mathrm{Ein}_{\rho}+\Scal_{\rho},
\end{split}
\]
where $R_{ijkl}=g(R(e_i,e_j)e_k,e_l)$ is the Riemannian curvature, $``\mathrm{Weyl}"$ the conformal Weyl tensor, $``\mathrm{Ein}"$ the Einstein tensor and $``\Scal"$ the scalar curvature. The coefficient $\pi_{\rho}(c_2)$ of $\Scal$ is an action of the second Casimir element $c_2=\sum_{1\le i,j\le n}e_{ij}e_{ji}$ on $W_{\rho}$, which acts on each irreducible bundle by a constant. Next we explain the first term  of \eqref{eqn:square1}, {\it the standard Laplacian} or {\it the Lichnerowicz Laplacian}. On the bundle $S_{\rho}$, we define the second order Laplace type operator by
\begin{equation}
\Delta_{\rho}:=\nabla^{\ast}\nabla+\frac{1}{2}R_{\rho}. \label{eqn:standard}
\end{equation}
This operator coincides with the Hodge Laplacian $\Delta=dd^{\ast}+d^{\ast}d$ on $S_{\rho}=\Lambda^k(T^{\ast}M)$ and $\Delta=D^2-\Scal/8$ on the spinor bundle. On an irreducible compact symmetric space $G/K$, it corresponds to the second Casimir operator (non-negative operator) for $G$. Remark that $\Delta_{\rho}$ is not necessarily non-negative operator on a compact Riemannian manifold in general. One of important observation for the standard Laplacian is that, when we have a decomposition $S_{\rho}=S_{\lambda_1}\oplus \cdots \oplus S_{\lambda_N}$ as associated vector bundle to $\SO(M)$ or $\Spin(M)$, the standard Laplacian $\Delta_{\rho}$ acts on each bundle diagonally (for more detail, \cite{SW2004}, \cite{H2}). 

We return to the case of $D(j)$. By the above observation, we have
\begin{equation*}
D(j)^2
=\begin{pmatrix}
\Delta_j+\frac{\Scal}{8} &  0 \\
  0     & \Delta_{j-1}+\frac{\Scal}{8} 
\end{pmatrix}-\frac{1}{2}
\begin{pmatrix}
(\id\otimes R_{\mathrm{Sym}^j_0})^j_j & (\id\otimes R_{\mathrm{Sym}^j_0})^{j-1}_{j} \\
(\id\otimes R_{\mathrm{Sym}^j_0})^{j}_{j-1} &(\id\otimes R_{\mathrm{Sym}^j_0})^{j-1}_{j-1}
\end{pmatrix}
\end{equation*}
along the decomposition $S_j\oplus S_{j-1}$. On the other hand, we take the square of the block matrix realization for $D(j)$, 
\begin{equation*}
D(j)^2=
\begin{pmatrix}
D_j^2+T_{j-1}^+T_j^- & D_jT_{j-1}^++T_{j-1}^+D_{j-1}'\\
T_j^-D_j+D_{j-1}'T_j^- &T_j^-T_{j-1}^++(D_{j-1}')^2
\end{pmatrix}. 
\end{equation*}
Comparing the above two equations, we show
\begin{gather*}
\begin{split}
D_j^2+T_{j-1}^+(T_{j-1}^+)^{\ast}=&\Delta_j+\frac{\Scal}{8}-\frac{1}{2}(\id\otimes R_{\mathrm{Sym}^j_0})^j_j, \\
 (T_{j-1}^+)^{\ast}T_{j-1}^++(D_{j-1}')^2=&\Delta_{j-1}+\frac{\Scal}{8}-\frac{1}{2}(\id\otimes R_{\mathrm{Sym}^j_0})^{j-1}_{j-1},
\end{split}  \\
(T_{j-1}^+)^{\ast}D_j+D_{j-1}'(T_{j-1}^+)^{\ast}=-\frac{1}{2}(\id\otimes R_{\mathrm{Sym}^j_0})^j_{j-1}, \quad
D_jT_{j-1}^++T_{j-1}^+D_{j-1}'=-\frac{1}{2}(\id\otimes R_{\mathrm{Sym}^j_0})^{j-1}_j. \end{gather*}
The curvature term $(\id\otimes R_{\mathrm{Sym}^j_0})^j_{j-1}$ is a bundle homomorphism from $S_{j}$ to $S_{j-1}$ depending only on the conformal Weyl tensor part $\mathrm{Weyl}_{\mathrm{Sym}^j_0}$ and the Einstein tensor part $\mathrm{Ein}_{\mathrm{Sym}^j_0}$ (\cite{H2}). Remark that, in the case of $j=1$, it depends only on the Einstein tensor, so that there are fruitful result of the Rarita-Schwinger operator on an Einstein manifold in \cite{HS}. On the other hand, the curvature term $(\id\otimes R_{\mathrm{Sym}^j_0})^j_j$ on $S_{j}$ depends not only on $\mathrm{Weyl}_{\mathrm{Sym}^j_0}$ and $\mathrm{Ein}_{\mathrm{Sym}^j_0}$ but also on the scalar curvature part. Therefore we can calculate the curvature terms explicitly on a Riemannian manifold of constant sectional curvature $K=c$ and know some operators commutes in a way. 
\begin{proposition}
If $(M,g)$ is a Riemannian manifold of constant sectional curvature $K=c$ with a spin structure, then
\begin{gather}
D_j^2+T_{j-1}^+(T_{j-1}^+)^{\ast}=\Delta_j-(j(n+j-2)-\frac{n(n-1)}{8})c, \label{eqn:j-1}\\
(T_{j-1}^+)^{\ast}T_{j-1}^++(D_{j-1}')^2=\Delta_{j-1}-(j(n+j-2)-\frac{n(n-1)}{8})c, \label{eqn:j-2}\\ 
(T_{j-1}^+)^{\ast}D_j+D_{j-1}'(T_{j-1}^+)^{\ast}=0,\quad D_jT_{j-1}^++T_{j-1}^+D_{j-1}'=0.\label{eqn:j-3}
\end{gather}
\end{proposition}
\begin{proof}
The conformal Weyl tensor and Einstein tensor are zero, and the scalar curvature is $\Scal=n(n-1)c$ on $(M,g)$. The action of $-\frac{1}{2}\id \otimes R_{\mathrm{Sym}^j_0}$ is constant by 
\[
-\frac{\pi_{\mathrm{Sym}^j_0}(c_2)}{2n(n-1)}\Scal=-j(n+j-2)c.
\]
Then we have
\[
\begin{split}
D_j^2+T_{j-1}^+T_j^-=\Delta_j-(j(n+j-2)-\frac{n(n-1)}{8})c.
\end{split}
\]
\end{proof}
Furthermore, taking $D(j)^2D(j)=D(j)D(j)^2$ into consideration, we conclude that
\[
\Delta_j D_j=D_j\Delta_j,\quad T_{j-1}^+\Delta_{j-1}=\Delta_{j}T_{j-1}^+,\quad T_j^-\Delta_j=\Delta_{j-1}T_j^-
\]
on a locally symmetric space $(M,g)$, whose curvature $R$ satisfies $\nabla R=0$. In particular, these commutation relations hold on a space of constant curvature. 

Now we calculate a normalizing constant $c_j$ in the equation $D_j=c_j\tilde{D}_j$. Thanks to Weitzenb\"ock formulas in \cite{H1} (see Appendix \ref{app:A} for more detail), the operators satisfy
\begin{align}
\frac{(n+2j)(n-2)}{n+2j-2}\tilde{D}_j^2+\frac{4(n+j-2)}{n+2j-2}(\tilde{T}_{j}^-)^{\ast}\tilde{T}_j^-&=\Delta_j+\mathrm{curv}, \label{WF1}\\
\frac{4j}{n+2j-2}(\tilde{T}_{j-1}^+)^{\ast}\tilde{T}_{j-1}^++\frac{(n+2j-4)(n-2)}{n+2j-2}(\tilde{D}_{j-1})^2&=\Delta_{j-1}+\mathrm{curv}, \label{WF2}
\end{align}
where ``curv'' means a bundle endomorphism depending on the Riemannian curvature as before. We also know that $D_j^2$ and $T_{j-1}^+T_j^-$ (resp. $\tilde{T}_{j+1}^-\tilde{T}_{j}^+$) are independent in the sense that there is no relation of the form
\[
a D_j^2+ b T_{j-1}^+T_j^-=\mathrm{curv},\quad (a,b)\neq (0,0).
\]
Then, comparing the equations \eqref{eqn:j-1},\eqref{eqn:j-2}, \eqref{WF1} and \eqref{WF2}, we obtain the normalizing constant for $D_j$ with respect to $\tilde{D}_j=\Pi_j\circ \nabla$. 
\begin{proposition}
We consider a Riemannian spin manifold $(M,g)$ and the vector bundle $S_j=\Spin(M)\times_{\Spin(n)} W_j$ with spin $j+1/2$. Let $\tilde{D}_j=\Pi_j\circ \nabla$ be the generalized gradient on $S_j$ defined by \eqref{h-Dirac}, and $D_j$ (resp. $D_{j-1}'$) be the operator defined by restricting the twisted Dirac operator $D(j)$ on $S_j$ (resp. on $S_{j-1}$) in \eqref{twisted-Dirac}. Then 
\[
D_j=\sqrt{\frac{(n+2j)(n-2)}{n+2j-2}}\tilde{D}_j,\quad D_{j-1}'=-\sqrt{\frac{(n+2(j-1)-2)(n-2)}{n+2(j-1)}}\tilde{D}_{j-1},
\]
and hence,
\[
D_j^2=\frac{(n+2j)(n-2)}{n+2j-2}\tilde{D}_j^2,\quad (D_{j-1}')^2=\frac{(n+2(j-1)-2)(n-2)}{n+2(j-1)}\tilde{D}_{j-1}^2. 
\]
As a result, we have
\[
D_j'=-\frac{n+2j-2}{n+2j}D_j.
\]
We also have normalizing constants for the twistor operators $T_j^{\pm}$ on $S_j$, 
\[
T_j^-=2\sqrt{\frac{n+j-2}{n+2j-2}}\tilde{T}_j^-,\quad T_j^+=2\sqrt{\frac{j+1}{n+2j}}\tilde{T}_j^+. 
\]
Hence we have
\begin{align*}
\nabla^{\ast}\nabla=&(\tilde{T}_j^+)^{\ast}\tilde{T}_j^++ U_j^{\ast}U_j+\tilde{D}_j^2+(\tilde{T}_j^-)^{\ast}\tilde{T}_j^-, \\
 =&\frac{n+2j}{4(j+1)}(T_j^+)^{\ast}T_j^++ U_j^{\ast}U_j+\frac{n+2j-2}{(n+2j)(n-2)}D_j^2+\frac{n+2j-2}{4(n+j-2)}(T_j^-)^{\ast}T_j^-. 
\end{align*}
\end{proposition}

\begin{rem}
The generalized gradient is defined up to constant multiple of a complex number $u$ with norm $|u|=1$, so that the above normalizing constant is unique up to such $u$. In the case of $D_j$ (resp. $D_j'$), taking into account the self adjointness, we put $u=1$ (resp. $u=-1$). 
\end{rem}
\begin{rem}
To calculate the normalizing constant for $D_{j}'=cD_j$, we may use the homomorphism type Weitzenb\"ock formula given in \cite{H2}. In fact, we have
\[
\frac{1}{\sqrt{n/2+j-2}}\tilde{D}_j\tilde{T}_{j-1}^+=\frac{1}{\sqrt{n/2+j}}\tilde{T}_{j-1}^+\tilde{D}_{j-1}.
\]
Therefore
\[
 D_jT_{j-1}^+=\frac{n+2j-4}{n+2j-2}T_{j-1}^+D_{j-1}.
\]
Then we have concluded $D_j'=-\frac{n+2j-2}{n+2j}D_j$. 
\end{rem}
\begin{exam}
The operator $D_0$ is just the Dirac operator $D$. The operator $D_1$ is the Rarita-Schwinger operator $Q$ in \cite{HS}.
\end{exam}
As a corollary, we rewrite the Weitzenb\"ock formula by $D_j$ and $T_j^{\pm}$.
\begin{corollary}\label{cor:W}
Let $(M,g)$ be a space of constant curvature $K=c$ with a spin structure. Then the following identities for operators from $\Gamma(S_j)$ to $\Gamma(S_j)$ hold:
\begin{equation}
\begin{split}
\Delta_j=&\nabla^{\ast}\nabla+\frac{1}{2}R_j=\nabla^{\ast}\nabla+(j(n+j-1)+\frac{n(n-1)}{8})c \\
  =&D_j^2+(T_j^-)^{\ast}T_j^- + (j(n+j-2)-\frac{n(n-1)}{8})c \\
  =&(T_{j}^+)^{\ast}T_{j}^++\frac{(n+2j-2)^2}{(n+2j)^2}D_{j}^2+((j+1)(n+j-1)-\frac{n(n-1)}{8})c. 
\end{split}\label{eqn:W-1}
\end{equation}
The operators among $\Gamma(S_{j-1})$ and $\Gamma(S_{j})$ are related as follows:
\begin{equation}
T_{j}^-D_j=\frac{n+2j-4}{n+2j-2}D_{j-1}T_j^-,\quad D_jT_{j-1}^+=\frac{n+2j-4}{n+2j-2}T_{j-1}^+D_{j-1}, \quad T_{j-1}^+=(T_{j}^-)^{\ast}.\label{TW}
\end{equation}
Moreover the operators commute with the standard Laplacian,
\[D_j\Delta_j=\Delta_jD_j, \quad T_j^-\Delta_j=\Delta_{j-1}T_j^-, \quad T_{j-1}^+\Delta_{j-1}=\Delta_{j}T_{j-1}^+.
\]
\end{corollary}
The main theorem in this section is a factorization formula for the power of the Laplacian on a space of constant curvature. On the flat space $\mathbb{R}^n$, it is known in \cite{BSSL1}, \cite{BSSL2} that $D_j$ factors through $(\Delta_j)^{j+1}$. In other words, there is a differential operator $A_{2j+1}$ with order $2j+1$ such that $(\Delta_j)^{j+1}=D_j\circ A_{2j+1}$. It follows that the known fundamental solution $G$ for $(\Delta_j)^{j+1}$ gives a fundamental solution $A_{2j+1}(G)$ for $D_j$. We can generalize this factorization on $\mathbb{R}^n$ to a space of constant curvature explicitly. 
\begin{theorem}[Factorization formula]\label{thm:fact}
Let $(M,g)$ be a space of constant curvature $K=c$ with a spin structure. For the higher spin Dirac operator $D_j$ and the standard Laplacian $\Delta_j$ on $S_j$, we define a second order operator on $\Gamma(S_j)$, 
\begin{equation*}
B(s;j):=D_j^2-\frac{(n+2s-2)^2}{(n+2j-2)^2}\left(\Delta_j-\left(s(n+s-2)-\frac{n(n-1)}{8}\right)c\right)
\end{equation*}
for $s=0,1,\cdots$. Then we have 
\begin{equation}
\prod_{s=0}^jB(s;j)=0. \label{eqn:fact1}
\end{equation}
\end{theorem}
\begin{proof}
We prove theorem by induction for $j$. We start from the equation for $j=0$,
\[
D_0^2-(\Delta_0+\frac{n(n-1)}{8}c)=0.
\]
This is just the Lichnerowicz formula for the Dirac operator. We assume the equation \eqref{eqn:fact1} holds for $j$. From \eqref{TW}, we easily show the key relations
\begin{equation}
\begin{split}
T_j^+B(s;j)=\frac{(n+2j)^2}{(n+2j-2)^2}B(s;j+1) T_j^+, &\quad 
T_{j+1}^-B(s;j+1)=\frac{(n+2j-2)^2}{(n+2j)^2}B(s;j)T_{j+1}^-,\\
(T_j^-)^{\ast}T_j^-=-B(j;j), \quad &\quad (T_j^+)T_j^+=-\frac{(n+2j-2)^2}{(n+2j)^2}B(j+1;j).
\end{split}
 \label{eqn:key}
\end{equation}
Sandwich \eqref{eqn:fact1} by $T_j^+$ and $(T_j^+)^{\ast}$, and we obtain
\[
\begin{split}
0=&T_j^+ \left(\prod_{s=0}^j B(s;j)\right)(T_j^+)^{\ast} =\left(\prod_{s=0}^j \frac{(n+2j)^2}{(n+2j-2)^2}B(s;j+1)\right) T_j^+(T_j^+)^{\ast},\\
=&-\frac{(n+2j)^{2(j+1)}}{(n+2j-2)^{2(j+1)}} \left(\prod_{s=0}^jB(s;j+1)\right)    B(j+1;j+1).
\end{split}
\]
Thus we have proved \eqref{eqn:fact1} holds for $j+1$. 
\end{proof}
To clarify the meaning of the above factorization formula, we shall introduce a filtration on $\Gamma(S_j)$. Put 
\[
F_j:=\ker T_j^-,\quad F_{j-1}:=\ker T_{j-1}^-T_j^-,\quad \cdots \quad F_1:=\ker T_1^-\cdots T_j^-,\quad F_0:=\Gamma(S_j),
\]
and we have a filtration on $\Gamma(S_j)$, 
\[
F_j\subset F_{j-1}\subset \cdots \subset F_1\subset F_0=\Gamma(S_j).
\]
\begin{proposition}\label{prop:6}On a compact space $(M,g)$ of constant curvature $K=c$ with a spin structure, 
\[
F_j=\ker B(j;j),\quad \cdots,\quad  F_k=\ker \prod_{s=k}^j B(s;j),\quad \cdots, \quad F_0=\Gamma(S_j). 
\]
\end{proposition}
\begin{proof}
If $\phi$ is in $F_k$, then $(T_{j}^-)^{\ast}\cdots (T_{k}^-)^{\ast}T_{k}^-\cdots T_j^-\phi=0$. Conversely, if $(T_{j}^-)^{\ast} \cdots (T_{k}^-)^{\ast}T_{k}^-\cdots T_j^-\phi$ is zero, then 
\[
\begin{split}
0=\int_M\langle (T_{j}^-)^{\ast}\cdots (T_{k}^-)^{\ast}T_{k}^-\cdots T_j^-\phi, \phi\rangle \mathrm{vol}_g =\int_M|T_{k}^-\cdots T_j^-\phi|^2\mathrm{vol}_g=\|T_{k}^-\cdots T_j^-\phi\|^2. 
\end{split}
\]
Thus $\phi$ is in $F_k$ if and only if $(T_{j}^-)^{\ast}\cdots (T_{k}^-)^{\ast}T_{k}^-\cdots T_j^-\phi=0$. From \eqref{eqn:key}, 
\[\begin{split}
(T_{j}^-)^{\ast} \cdots (T_{k}^-)^{\ast}T_{k}^- \cdots T_j^-=&-(T_{j}^-)^{\ast} \cdots (T_{k+1}^-)^{\ast}B(k;k) T_{k+1}^-\cdots T_j^-\\
=&cB(k;j)(T_{j}^-)^{\ast} \cdots (T_{k+1}^-)^{\ast}T_{k+1}^-\cdots T_j^- \quad (\exists c\neq 0)\\
=&-cB(k;j)(T_{j}^-)^{\ast} \cdots (T_{k+2}^-)^{\ast}B(k+1;k+1)T_{k+2}^-\cdots T_j^-\\
  &\cdots \\
=&c'B(k;j)\cdots B(j;j),\quad (\exists c'\neq 0)
\end{split}
\]
Then we can prove the proposition. 
\end{proof}
The above natural filtration allows us to get an associated grading on $\Gamma(S_j)$. In other words, there are subspaces $W_s$ for $s=0,\cdots,j$ such that $F_s=W_s\oplus F_{s+1}$ and $\Gamma(S_j)=\oplus_{0\le s\le j}W_s$. 
\begin{theorem}\label{thm:7}
Let $S_j$ be the bundle for spinor fields with spin $j+1/2$ on a compact space of constant curvature $K=c$ with a spin structure. We put
\[
W_s:=\begin{cases} 
 T_{j-1}^+\cdots T_0^+(\Gamma(S_0)), & s=0,\\
T_{j-1}^+\cdots T_{s}^+(\ker T_{s}^-), & 1\le s\le j-1, \\
 \ker T_j^-,  & s=j.
\end{cases}
\]
Then the space $\Gamma(S_j)$ has an orthogonal decomposition associated to the filtration $\{F_j\}_j$, 
\[
\Gamma(S_j)=\bigoplus_{0\le s\le j}W_s,\quad F_s=W_s\oplus F_{s+1} \quad (0\le s\le j).
\]
Furthermore, $W_s\subset \ker B(s;j)$ for $s=0,\cdots,j$. 
\end{theorem}
\begin{proof}

From Proposition \ref{prop:1}, there are the orthogonal decompositions
\[
\Gamma(S_s)=\ker T_s^-\oplus T_{s-1}^+(\Gamma(S_{s-1}))
\]
for $s=0,1,\cdots$. When we act $T_s^+$ on $\Gamma(S_s)$, we can prove the decomposition is preserved, 
\[
T_s^+(\Gamma(S_s))=T_s^+(\ker T_s^-)\oplus T_s^+T_{s-1}^+(\Gamma(S_{s-1})).
\]

Since $T_s^+(\Gamma(S_s))=T_s^+(\ker T_s^-)+T_s^+T_{s-1}^+(\Gamma(S_{s-1}))$, we only have to prove the components intersect orthogonally. Take $\phi_s$ in $\ker T_s^-$, and we show 
\[
\begin{split}
 T_s^-T_{s+1}^-T_s^+(\phi_s)=&T_s^-(T_s^+)^{\ast}T_s^+\phi_s =-\frac{(n+2s-2)^2}{(n+2s)^2}T_s^-B(s+1;s)\phi_s\\
=&-\frac{(n+2s-4)^2}{(n+2s)^2}B(s+1;s-1)T_s^-\phi_s=0.
\end{split}
\]
Then 
\[
(T_s^+\phi_s,T_s^+T_{s-1}^+\psi_{s-1})=(T_s^-T_{s+1}^-T_s^+\phi_s,\psi_{s-1})=0 \quad \textrm{for $\psi_{s-1}\in \Gamma(S_{s-1})$}. 
\]
Thus we have obtained an orthogonal decomposition $T_s^+(\Gamma(S_s))=T_s^+(\ker T_s^-)\oplus T_s^+T_{s-1}^+(\Gamma(S_{s-1}))$. As a consequence, we have
\[
\begin{split}
\Gamma(S_{s+1})=\ker T_{s+1}^-\oplus T_s^+(\ker T_s^-)\oplus T_s^+T_{s-1}^+(\Gamma(S_{s-1})), \\
F_{s+1}=\ker T_{s+1}^-,\quad F_{s}=\ker T_{s+1}^-\oplus T_s^+(\ker T_s^-).
\end{split}
\]
Repeating this, we have concluded $\Gamma(S_j)=\oplus_s W_s$. 

Next we shall prove $W_s\subset \ker B(s;j)$. For $\phi_s$ in $\ker T_s^-$, 
\[
\begin{split}
B(s;j)T_{j-1}^+\cdots T_s^+\phi_s=&\frac{(n+2j-4)^2}{(n+2j-2)^2}T_{j-1}^+B(s;j-1)T_{j-2}^+ \cdots T_s^+\phi_s\\
=&\frac{(n+2s-2)^2}{(n+2j-2)^2}T_{j-1}^+\cdots T_{s+1}^+T_s^+B(s;s)\phi_s\\
=&-\frac{(n+2s-2)^2}{(n+2j-2)^2}T_{j-1}^+\cdots T_{s+1}^+T_s^+(T_s^-)^{\ast}T_s^-\phi_s=0.
\end{split}
\]
Thus we have concluded that $W_s\subset \ker B(s;j)$. Remark that the reverse inclusion doesn't hold in general (for example, the case of $K=0$). 
\end{proof}

Since $D_j\Delta_j=\Delta_jD_j$, the operator $D_j$ and $\Delta_j$ have simultaneous eigenstates on a compact space $(M,g)$ of constant curvature and their eigenvalues relate to each other. Let $\phi$ is a simultaneous eigenspinor in $W_{k}$ with eigenvalue $\mu$ for $\Delta_j$ and $\lambda^2$ for $D_j^2$. Then 
\[
\lambda^2=\frac{(n+2k-2)^2}{(n+2j-2)^2}\left(\mu-\left(k(n+k-2)-\frac{n(n-1)}{8}\right)c\right).
\]
For the non-compact case, the situation is more difficult. In fact, it has already been known that $\Delta_j\phi=0$ doesn't follow $D_j\phi=0$ on $\mathbb{R}^n$, and vice versa (\cite{BSSL2}).

At the end of this section, we state some remarks.
\begin{rem}
It was shown in \cite{HS} that the Rarita-Schwinger operator $D_1$ satisfies
\[
\left(D_1^2-\left(\Delta_1+\frac{n-8}{8n}\Scal \right)\right)\left(D_1^2-
\frac{(n-2)^2}{n^2}\left(\Delta_1+\frac{\Scal}{8}\right)\right)=0
\]
on an Einstein manifold $(M,g)$. To get a factorization formula for higher spin cases of $j>1$, we have to assume the sectional curvature is constant because the curvature term and their derivatives remain in Weitzenb\"ock formulas and the formula is too complicated to get a factorization. To calculate such error terms is a future problem. 
\end{rem}
\begin{rem}\label{symbol}
The above  theorem  shows that $\sigma_{\xi}(D_j)$ is a factor of $\sigma_{\xi}(\Delta_j^{j+1})=(-\|\xi\|^2)^{j+1}$ on the principal symbol level, and so that $D_j$ is elliptic. 
\end{rem}
\begin{rem}
We can rewrite the equation \eqref{eqn:fact1} for $T_j^+$ as
\[
\prod_{s=0}^j\left((T_j^+)^{\ast}T_j^+-a'(s;j)(\Delta_j-b'(s;j)c)\right)=0,
\]
where
\[
a'(s;j)=\frac{4(j-s+1)(n+j+s-1)}{(n+2j)^2},\quad b'(s;j)=j(n+j)+s(n+s-2)+\frac{n(n+1)}{8}.
\]
Then the eigenvalue of the standard Laplacian $\Delta_j$ restricted to $\ker T_j^+\cap W_k$ on a compact space of constant curvature is just
\[
(j(n+j)+k(n+k-2)+\frac{n(n+1)}{8})c.
\]
On the other hand, the formula for $T_j^-$ is
\[
\left(\prod_{s=1}^j\left((T_j^-)^{\ast}T_j^--a''(s;j)(\Delta_j-b''(s;j)c)\right)\right)(T_j^-)^{\ast}T_j^-=0.
\]
In this case, we can not say the similar statement for $\phi$ in $\ker T_j^-$. 
\end{rem}

\section{Harmonic analysis for spinor fields with higher spin on $S^n$}\label{sec:3}

In this section, we study harmonic analysis for spin $j+1/2$ fields on the sphere. We calculate the eigenvalues of $D_j^2$ and other related second order operators and show how the space of spinor fields are related to each other as $\Spin(n+1)$-modules.
 
Let $S^n$ be the $n$-dimensional standard sphere. We realize it as a $\Spin(n+1)$-orbit of $e_{n+1}=(0,\cdots,0,1)$ in $\mathbb{R}^{n+1}$. The isotropy group is 
\[
\Spin(n)=\{g\in \Spin(n+1)|g\cdot e_{n+1}=e_{n+1}\},
\]
and the sphere is a symmetric space $\Spin(n+1)/\Spin(n)$. The (unique) spin structure of $S^n$ is the principal $\Spin(n)$-bundle $\Spin(n+1)$ where $\Spin(n)$ acts naturally from the right side. The bundle for higher spinor fields is given by
\[
S_j=\Spin(n+1)\times_{\Spin(n)}W_j
\]
with a natural fiber metric induced from a $\Spin(n)$-invariant inner product on $W_j$. The space $L^2(S^n,S_j)$ of $L^2$-integrable sections of $S_j$ is a unitary representation space of $\Spin(n+1)$ with respect to $L^2$-inner product $(\phi,\psi)=\int_{S^n}\langle \phi,\psi\rangle \mathrm{vol}_g$. Frobenius reciprocity provides its irreducible decomposition, 
\[
L^2(S^n,S_j)=\bigoplus_{\rho \in \widehat{\Spin(n+1)}} V_{\rho}\otimes \Hom_{\Spin(n)}(V_{\rho},W_j),
\]
where $\widehat{\Spin(n+1)}$ means the equivalent classes of irreducible unitary finite dimensional representations of $\Spin(n+1)$. Such class is parametrized by its highest weight $\rho=(\rho^1,\dots,\rho^m)$ for $m=[\frac{n+1}{2}]$ satisfying $\rho$ in $\mathbb{Z}^m\cup (\mathbb{Z}+1/2)^m$ and the dominant condition
\begin{align*}
\rho^1\ge \rho^2\ge \cdots \ge \rho^m\ge 0, \quad (\textrm{for $n+1=2m+1$}),\\
\rho^1\ge \rho^2\ge \cdots \ge \rho^{m-1}\ge |\rho^m|, \quad (\textrm{for $n+1=2m$}).
\end{align*}
We state the branching rule to describe the restriction of an irreducible representation of $\Spin(n+1)$ to $\Spin(n)$ (cf. \cite{Z}). Fix an irreducible $\Spin(n+1)$-module $V_{\rho}$ with the highest weight $\rho$. When $V_{\rho}$ is considered as a $\Spin(n)$-module, the following $\Spin(n)$-modules appear with each multiplicity $1$ as irreducible summands:
\begin{itemize}
\item In the case of $n+1=2m+1$, $\lambda=(\lambda^1,\cdots,\lambda^{m})$ in $\widehat{\Spin(2m)}$ satisfying
\[
\rho^1\ge \lambda^1 \ge \rho^2\ge \lambda^2\ge \cdots \ge \rho^{m}\ge \lambda^m\ge -\rho^m, \qquad  \rho^1-\lambda^1\in \mathbb{Z}.
\]
\item In the case of $n+1=2m$,  $\lambda=(\lambda^1,\cdots,\lambda^{m-1})$ in $\widehat{\Spin(2m-1)}$ satisfying
\[
\rho^1\ge \lambda^1 \ge \rho^2\ge \lambda^2\ge \cdots \ge \rho^{m-1}\ge \lambda^{m-1}\ge |\rho^m|, \qquad \rho^1-\lambda^1\in \mathbb{Z}.
\]
\end{itemize}
We denote by $V_j(k,s)$ an irreducible $\SO(n+1)$-module with the highest weight
\[
(k+j,s,0_{m-2})
\]
for $k$ in $\mathbb{Z}_{\ge 0}$ and $s=0,\dots, j$. We also denote by $V_j(k,s)'$ an irreducible $\Spin(2m+1)$-module with the highest weight
\[(k+j+1/2,s+1/2,(1/2)_{m-2}),
\]
and by $V_j^{\pm}(k,s)'$ an irreducible $\Spin(2m)$-module with the highest weight 
\[
(k+j+1/2,s+1/2,(1/2)_{m-3},\pm 1/2)
\]
for $k$ in $\mathbb{Z}_{\ge 0}$ and $s=0,\dots, j$. 
\begin{rem}
For the case of $n+1=4$, the space $V_j(k,s)$ splits the sum of $\SO(4)$-modules  with the highest weight $(k+j,s)$ and $(k+j,-s)$. The space $V_j^{\pm}(k,s)'$ means an irreducible $\Spin(4)$-module with the highest weight $(k+j+1/2,\pm (s+1/2))$.
\end{rem}

Applying Frobenius reciprocity and the branching rule to our case $\lambda=(j+1/2,1/2,\cdots,1/2)$, we have the proposition.
\begin{proposition}
\begin{enumerate}
\item  For $n=2m$,
\begin{align*}
L^2(S^{2m},S_j^{\pm})=\bigoplus_{0\le s\le j} \bigoplus_{k\in \mathbb{Z}_{\ge 0}}V_j(k,s)',\quad L^2(S^{2m},S_j)=\bigoplus_{0\le s\le j} \bigoplus_{k\in \mathbb{Z}_{\ge 0}} 2V_j(k,s)'.
\end{align*}
\item For $n=2m-1$,
\[
L^2(S^{2m-1},S_j)=\bigoplus_{0\le s\le j} \bigoplus_{k\in \mathbb{Z}_{\ge 0}} V_j^+(k,s)'\oplus V_j^-(k,s)'.
\]
\end{enumerate}
\end{proposition}
From now on, we put $V_j(k,s)'=V_j^+(k,s)'\oplus V_j^-(k,s)'$ for $n=2m-1$. By Weyl's dimension formula, we have
\[
\dim V_j(k,s)'=2^{[\frac{n+1}{2}]}\frac{(l+n-1+s)(l+1-s)}{(n-1)(n-2)}\binom{l+n-2}{l+1}\binom{s+n-3}{s},
\]
where we set $l=k+j$. To calculate the eigenvalues of the square of the higher spin Dirac operator $D_j^2$, we need the ones of $\Delta_j$ with adding a constant. 
\begin{lemma}\label{lem-delta}
On $V_j(k,s)'$, it holds that
\[
\left.\left(\Delta_j-(s(s+n-2)-\frac{n(n-1)}{8})\right)\right|_{V_j(k,s)'}=\left(j+k+\frac{n}{2}\right)^2.
\]
\end{lemma}
\begin{proof}
The standard Laplacian $\Delta_j$ on $S_j$ coincides with the Casimir operator $c_2/2$ for $\Spin(n+1)$. By Freudenthal's formula, the eigenvalue of the $c_2/2$ is $\pi_{\rho}(c_2/2)=\la \rho,\rho\ra+2\la \rho,\delta_{\Spin(n+1)}\ra$. The inner product is given by $\la \mu,\nu\ra=\sum_{1\le i\le m} \mu^i\nu^i$ and $\delta_{\Spin(n+1)}$ is half the sum of the positive roots,
\[
\delta_{\Spin(n+1)}=
\begin{cases}
(m-1/2,m-2/3,\cdots,3/2,1/2)     &     \textrm{for $n+1=2m+1$},\\
(m-1,m-2,\cdots, 1,0)   & \textrm{for $n+1=2m$}.
\end{cases}
\]
Then
\[
\begin{split}
 &\pi_{\rho}(c_2/2)-(s(s+n-2)-\frac{n(n-1)}{8})\\
=&(k+j+1/2)(k+j+n-1/2)+(s+1/2)(s+n-5/2)\\
  &\quad +(n-4)(n-3)/8-(s(s+n-2)-\frac{n(n-1)}{8})\\
=&(j+k+n/2)^2.
\end{split}
\]
\end{proof}
From the factorization formula in Theorem \ref{thm:fact}, we can guess the eigenvalue $\lambda^2$ of $D_j^2$ on $V_j(k,s)'$ is $\frac{(n+2s-2)^2}{(n+2j-2)^2}\left(j+k+\frac{n}{2}\right)^2$. In fact, we have 
\begin{theorem}
On the $n$-dimensional standard sphere $S^n$, the eigenvalue of $D_j^2$ on $V_j(k,s)'$ is 
\[
\frac{(n+2s-2)^2}{(n+2j-2)^2}\left(j+k+\frac{n}{2}\right)^2
\]
with multiplicity $\dim V_j(k,s)'$. 
\end{theorem}
\begin{proof}
We shall prove the case of $n=2m-1$ by induction for $j\ge 0$. First, we consider the case of $j=0$. It is well-known that the spectrum of the Dirac operator $D_0^2$ on the sphere is $(k+n/2)^2$ on $V_0(k,0)'$ for $k=0,1,2,\dots$. Next, assuming the statement in theorem for $j$ holds, we shall prove the one for $j+1$. We take $\phi$ in $V_j(k,s)'$. It follows from our assumption, Lemma \ref{lem-delta} and \eqref{eqn:W-1} that 
\[
(T_j^+)^{\ast}T_j^+\phi=\frac{k (j-s+1) (2 j+k+n) (j+s+n-1)}{(j+n/2)^2}\phi.
\]
Then 
\[
\ker T_j^+=\bigoplus_{0\le s\le j} V_j(0,s)',
\]
and $T_j^+:V_j(k,s)'\to V_{j+1}(k-1,s)'$ for any $k\ge 1$ is an isomorphism as $\Spin(n+1)$-module. Here we use Schur's lemma with respect to $\Spin(n+1)$-invariant operator $T_j^+$. Taking such a $T_j^+\phi$ in $V_{j+1}(k-1,s)'$, we show
\[
D_{j+1}^2T_j^+\phi=\frac{(n+2j-2)^2}{(n+2j)^2}T_j^+D_j^2\phi=\frac{(n+2s-2)^2}{(n+2j)^2}\left(j+1+(k-1)+\frac{n}{2}\right)^2T_j^+\phi.
\]
The case of $s=j+1$ remains to be proved. According to Proposition \ref{prop:1}, $L^2(S_{j+1})$ decomposes as $\ker T_{j+1}^-\oplus \mathrm{Image}\;T_j^+$. Then 
\[
\ker T_{j+1}^-=\bigoplus_{k\ge 0} V_{j+1}(k,j+1)'.
\]
When $\phi$ is in $V_{j+1}(k,j+1)'$, Weitzenb\"ock formula \eqref{eqn:W-1} gives
\[
D_{j+1}^2\phi=\Delta_{j+1}\phi-((j+1)(n+j-1)-\frac{n(n-1)}{8})\phi=\left((j+1)+k+\frac{n}{2}\right)^2\phi.
\]
Thus we have proved the statement for $j+1$. We can prove the case of $n=2m$ when we double each component $V_j(k,s)'$ in the above proof.
\end{proof}
By using Corollary \ref{cor:W}, we calculate all the eigenvalues of the other generalized gradients on $S_j$. 
\begin{corollary}
The eigenvalues of $(T_j^+)^{\ast}T_j^+$, $(T_j^-)^{\ast}T_j^-$ and $U_j^{\ast}U_j$ on $V_j(k,s)'$ are given by 
\begin{gather}
\frac{k (j-s+1) (2 j+k+n) (j+s+n-1)}{(j+n/2)^2},\\
 \frac{(k+1) (j-s) (2j+k+n-1) (j+s+n-2)}{(j+n/2-1)^2},\\
\textrm{and}\quad \frac{(n-3) s (j+k+1) (n+s-2) (j+k+n-1)}{(n-2)(j+1)(j+n-2)},
\end{gather}
respectively. 

\end{corollary}
We investigate how the eigenspaces $\{V_j(k,s)'\}_{k,j,s}$ relate to each other in virtue of the generalized gradients. Put 
\[
\textbf{V}_j(s)':=\bigoplus_{k\in \mathbb{Z}_{\ge 0}}V_j(k,s)'
\]
for each $j$ and $s=0,\cdots,j$, and we show that $\textbf{V}_j(s)'$ is just $W_s$ in Theorem \ref{thm:7}. Remark that, for the case of the sphere, $W_s$ coincides with the kernel of $B(s;j)$. The following diagram is useful to understand relations between the spaces: 
\[
\begin{split}
&\xymatrix@C=3pt{
L^2(S_0) \quad= & \mathbf{V}_0(0)' \ar@{->>}[d]^{T_0^+} &          && &  && \\
          & \vdots  & \ddots    && && &}\\
&\xymatrix@C=3pt{
L^2(S_{j-2})=  & \mathbf{V}_{j-2}(0)'  \ar@<1ex>@{->>}[d]^{T_{j-2}^+}  & \oplus \cdots \oplus    & \mathbf{V}_{j-2}(j-2)' \ar@<1ex>@{->>}[d]^{T_{j-2}^+}    & \oplus  & \{0\} &&   \\
L^2(S_{j-1})=  & \mathbf{V}_{j-1}(0)' \ar@{^{(}-_>}[u]^{T_{j-1}^-} \ar@<1ex>@{->>}[d]^{T_{j-1}^+} & \oplus \cdots \oplus   & \mathbf{V}_{j-1}(j-2)' \ar@{^{(}-_>}[u]^{T_{j-1}^-} \ar@<1ex>@{->>}[d]^{T_{j-1}^+} & \oplus & \mathbf{V}_{j-1}(j-1)' \ar[u]^{T_{j-1}^-} \ar@<1ex>@{->>}[d]^{T_{j-1}^+}&  \oplus & \{0\}   \\
L^2(S_j)= & \mathbf{V}_{j}(0)' \ar@{^{(}-_>}[u]^{T_{j}^-} \ar@<1ex>@{->>}[d]^{T_{j}^+}  & \oplus \cdots \oplus   & \mathbf{V}_{j}(j-2)' \ar@{^{(}-_>}[u]^{T_{j}^-} \ar@<1ex>@{->>}[d]^{T_{j}^+}& \oplus &
 \mathbf{V}_{j}(j-1)' \ar@{^{(}-_>}[u]^{T_{j}^-} \ar@<1ex>@{->>}[d]^{T_{j}^+} & \oplus & \mathbf{V}_{j}(j)' \ar[u]^{T_{j}^-}\ar@<1ex>@{->>}[d]^{T_{j}^+}\\
& \ar@{^{(}-_>}[u]^{T_{j+1}^-}& & \ar@{^{(}-_>}[u]^{T_{j+1}^-}& &\ar@{^{(}-_>}[u]^{T_{j+1}^-}&&\ar@{^{(}-_>}[u]^{T_{j+1}^-}
}
\end{split}
\]
Here we change each $\mathbf{V}_j(s)'$ by $2\mathbf{V}_j(s)'$ when $n$ is $2m$. Then we have the proposition. 
\begin{proposition}
The kernels and images of the generalized gradients are realized as $\Spin(n+1)$-module in $L^2(S^n,S_j)$ for $j=0,1,\cdots$ as follows. 

For $n=2m-1$, 
\begin{gather*}
 \ker T_j^+  =\bigoplus_{0\le s\le j} V_j(0,s)', \quad 
\ker T_j^-= \mathbf{V}_j(j)', \quad \ker U_j =\mathbf{V}_j(0)',\\
\mathrm{Image}\; T_{j-1}^+=\bigoplus_{0\le s\le j-1} \mathbf{V}_{j}(s)', \quad \mathrm{Image}\; T_{j+1}^-=L^2(S^n,S_j)\ominus \ker T_j^+.
\end{gather*}
In particular, we have
\[
\ker T_j^+\cap \ker T_j^-=V_j(0,j)', \qquad \ker T_j^+\cap \ker U_j=V_j(0,0)'.
\]
For $n=2m$, we double each component $V_j(k,s)'$ and $\mathbf{V}_j(s)'$ on the above equations.

\end{proposition}
\section{Generalized gradients on $j$-th symmetric tensor fields}\label{sec:4}
In this section, we study the generalized gradients on the bundle of trace-free symmetric tensors $\Sym^j_0=\Sym_0^j(TM^c)$. Let $(M,g)$ be an oriented Riemannian manifold. The bundle $\Sym^j(TM^c)$ for $j$-th symmetric tensor fields splits as \eqref{symj} and the primitive component $\Sym^j_0$ is a vector bundle associated to the (oriented) orthonormal frame bundle $\SO(M)$ with the highest weight $(j,0_{m-1})$. We compose the covariant derivative $\nabla$ and the orthogonal projection along the decomposition 
\[
\Sym_0^{j}\otimes TM^c=\Sym_0^{j,1}\oplus \Sym_0^{j+1}\oplus \Sym_0^{j-1},
\]
where $\Sym_0^{j,1}$ is an irreducible vector bundle with the highest weight $(j,1,0_{m-2})$. Then we have three generalized gradients 
\begin{align*}
T_j^+:\Gamma(\Sym_0^j)\to \Gamma(\Sym_0^{j+1}),\quad U_j:\Gamma(\Sym_0^j)\to \Gamma(\Sym_0^{j,1}),\quad T_j^-:\Gamma(\Sym_0^j)\to \Gamma(\Sym_0^{j-1}),
\end{align*}
where we set $U_0=0$ and $T_0^-=0$. 
\begin{rem}
For $n=4$, $\Sym_0^{j.1}$ splits into the sum of $\SO(4)$-modules with the highest weight $(j,1)$ and $(j,-1)$. For $n=3$ and $j\ge 1$, $\Sym_0^j$ appears again instead of $\Sym_0^{j,1}$. We denote by $U_j$ the self-gradient from $\Gamma(\Sym_0^j)$ to $\Gamma(\Sym_0^j)$ on a $3$-dim Riemannian manifold, which is not an elliptic operator.
\end{rem}
It follows from Weitzenb\"ock formula \cite{H1} that there are two identities among them
\begin{align*}
\nabla^{\ast}\nabla&=(T_j^+)^{\ast}T_j^++U_j^{\ast}U_j+(T_j^-)^{\ast}T_j^-,\\
\frac{1}{2}R_{\Sym^j_0}&=-j(T_j^+)^{\ast}T_j^++U_j^{\ast}U_j+(n+j-2)(T_j^-)^{\ast}T_j^-.
\end{align*}
Deleting $U_j$ from the above equations, we find
\[
\Delta_j:=\nabla^{\ast}\nabla+\frac{1}{2}R_{\Sym^j_0}=(j+1)(T_j^+)^{\ast}T_j^+-(n+j-3)(T_j^-)^{\ast}T_j^-+R_{\Sym^j_0}.
\]
As the case of the higher spin fields, we have to compare $T_{j+1}^-$ and $(T_j^+)^{\ast}$ from $\Gamma(\Sym_0^{j+1})$ to $\Gamma(\Sym_0^{j})$. The next lemma follows from \cite{H1}.
\begin{lemma}\label{relative}
Let $D^{\rho}_{\lambda}=\Pi_{\lambda}\circ \nabla :\Gamma(S_{\rho})\to \Gamma(S_{\lambda})$ be a generalized gradient on an irreducible vector bundle $S_{\rho}=\SO(M)\times_{SO(n)}W_{\rho}$ over a Riemannian manifold $(M,g)$. Here $\Pi_{\lambda}$ is the orthogonal projection onto $S_{\lambda}$ from $S_{\rho}\otimes TM^c$. Then 
\begin{equation}
(D^{\rho}_{\lambda})^{\ast}D^{\rho}_{\lambda}=\frac{\dim W_{\lambda}}{\dim W_{\rho}}D_{\rho}^{\lambda}(D_{\rho}^{\lambda})^{\ast}.\label{eqn:rel}
\end{equation}
\end{lemma}
\begin{proof}
The principal symbol of $D^{\rho}_{\lambda}$ is given by the linear map $p^{\rho}_{\lambda}(\xi)$ from $W_{\rho}$ to $W_{\lambda}$ defined by
\[
p^{\rho}_{\lambda}(\xi)\phi:=\Pi_{\lambda}(\phi\otimes \xi),\quad  \textrm{for $\phi$ in $W_{\rho}$ and $\xi$ in $\mathbb{C}^n= T_x^{\ast}M^c$}.
\]
In other words, $D^{\rho}_{\lambda}$ is $\sum_i p^{\rho}_{\lambda}(e_i)\nabla_{e_i}$ as the Dirac operator $D$ is given by $\sum_i e_i\cdot \nabla_{e_i}$. We denote by $p^{\rho}_{\lambda}(\xi)^{\ast}$ its adjoint map from  $V_{\lambda}$ to $V_{\rho}$.  Conversely we know that the orthogonal projection $\Pi_{\lambda}$ is realized by using the principal symbol (see Lemma 4.13 in \cite{H1}),
\[
\Pi_{\lambda}(\phi\otimes \xi)=\sum_{1\le i\le n} (p^{\rho}_{\lambda}(e_i)^{\ast}p^{\rho}_{\lambda}(\xi)\phi) \otimes e_i,
\]
where $\{e_i\}_i$ is an orthonormal basis for $\mathbb{C}^n$. We shall calculate the following two $\SO(n)$-invariant maps from $W_{\rho}$ to $W_{\rho}$, which are constants by Schur's lemma,
\[
\sum_i p^{\rho}_{\lambda}(e_i)^{\ast}p^{\rho}_{\lambda}(e_i),\quad \sum_i p^{\lambda}_{\rho}(e_i)p^{\lambda}_{\rho}(e_i)^{\ast}.
\]
Let $\{\phi_{\alpha}\}_{\alpha}$ be a basis of $W_{\rho}$. We consider the trace of $\Pi_{\lambda}:W_{\rho}\otimes \mathbb{C}^n\to W_{\lambda}$. 
\[
\begin{split}
\dim W_{\lambda}=&\sum_{\alpha,i} \la \Pi_{\lambda}(\phi_{\alpha}\otimes e_i),\phi_{\alpha}\otimes e_i\ra
   =\sum_{\alpha,i,j} \la (p^{\rho}_{\lambda}(e_j)^{\ast}p^{\rho}_{\lambda}(e_i)\phi_{\alpha}) \otimes e_j,\phi_{\alpha}\otimes e_i \ra\\
=&\sum_{\alpha,i} \la p^{\rho}_{\lambda}(e_i)^{\ast}p^{\rho}_{\lambda}(e_i)\phi_{\alpha}, \phi_{\alpha}\ra=\sum_i p^{\rho}_{\lambda}(e_i)^{\ast}p^{\rho}_{\lambda}(e_i) \sum_{\alpha} \la \phi_{\alpha}, \phi_{\alpha}\ra=\dim W_{\rho} \sum_i p^{\rho}_{\lambda}(e_i)^{\ast}p^{\rho}_{\lambda}(e_i). 
\end{split}
\]
Then we have $\sum_i p^{\rho}_{\lambda}(e_i)^{\ast}p^{\rho}_{\lambda}(e_i)=\frac{\dim W_{\lambda}}{\dim W_{\rho}}\id_{W_{\rho}}$. On the other hand, we have
\[
\begin{split}
 \la p^{\lambda}_{\rho}(\xi)\phi, p^{\lambda}_{\rho}(\eta)\psi\ra =&\la \Pi_{\rho}(\phi\otimes \xi),\Pi_{\rho}(\psi \otimes \eta)\ra
=\sum_{i,j} \la p^{\lambda}_{\rho}(e_i)^{\ast} p^{\lambda}_{\rho}(\xi)\phi\otimes e_i , p^{\lambda}_{\rho}(e_j)^{\ast}p^{\lambda}_{\rho}(\eta)\psi\otimes e_j\ra \\
=&\sum_{i} \la p^{\lambda}_{\rho}(e_i)^{\ast} p^{\lambda}_{\rho}(\xi)\phi, p^{\lambda}_{\rho}(e_i)^{\ast}p^{\lambda}_{\rho}(\eta)\psi\ra=\sum_{i} \la  p^{\lambda}_{\rho}(e_i) p^{\lambda}_{\rho}(e_i)^{\ast} p^{\lambda}_{\rho}(\xi)\phi, p^{\lambda}_{\rho}(\eta)\psi\ra\\
=&(\sum_{i} p^{\lambda}_{\rho}(e_i) p^{\lambda}_{\rho}(e_i)^{\ast})\la p^{\lambda}_{\rho}(\xi)\phi, p^{\lambda}_{\rho}(\eta)\psi\ra.
\end{split}
\]
Then we have $\sum_{i} p^{\lambda}_{\rho}(e_i) p^{\lambda}_{\rho}(e_i)^{\ast}=\id_{W_{\rho}}$. Since there is a nonzero constant $a$ such that $p^{\rho}_{\lambda}(\xi)=ap^{\lambda}_{\rho}(\xi)^{\ast}$, we have concluded that $|a|^2=\dim W_{\lambda}/\dim W_{\rho}$ and hence $(D^{\rho}_{\lambda})^{\ast}D^{\rho}_{\lambda}=\frac{\dim W_{\lambda}}{\dim W_{\rho}}D_{\rho}^{\lambda}(D_{\rho}^{\lambda})^{\ast}$.
\end{proof}

Thus we have given relations among the operators to construct a factorization formula on a space of constant curvature. 
\begin{proposition}\label{prop:W-sym}
On a space $(M,g)$ of constant curvature $K=c$, the operators $T_j^+$, $T_j^-$ and $\Delta_j$ for $j=0,1,2,\dots$ satisfy 
\begin{equation}
\begin{split}
\Delta_j&=(j+1)(T_j^+)^{\ast}T_j^+-(n+j-3)(T_j^-)^{\ast}T_j^-+2j(n+j-2)c,\\
\frac{1}{2}R_{\Sym^j_0}&=j(n+j-2)c=-j(T_j^+)^{\ast}T_j^++U_j^{\ast}U_j^++(n+j-2)(T_j^-)^{\ast}T_j^-, \\
(T_{j+1}^-)^{\ast}T_{j+1}^-&=\frac{(j+1)(n+2j-2)}{(n+j-2)(n+2j)}T_j^+(T_j^+)^{\ast}.
\end{split} \label{eqn:W-T}
\end{equation}
\end{proposition}
\begin{rem}
We don't normalize constant multiples of the operators as the previous section. When we use a notation in \cite{HMS}, we find
\[
(T_j^+)^{\ast}T_j^+=\frac{1}{j+1}d_0^{\ast}d_0,\quad (T_j^-)^{\ast}T_j^-=\frac{n+2j-4}{(n+2j-2)(n+j-3)}\delta^{\ast}\delta.
\]
\end{rem}
We show the factorization formula for the $j+1$-st power of the Laplacian on $j$-th symmetric tensor fields, which will be useful to calculate the eigenvalues on the sphere in the next section.  
\begin{theorem}[Factorization formula]\label{thm:fact2}
Let $(M,g)$ be a space of constant curvature $K=c$, and $\Sym_0^j$ the vector bundle for trace-free symmetric tensor fields on $(M,g)$. The operator $T_j^+:\Gamma(\Sym_0^j)\to \Gamma(\Sym_0^{j+1})$ and the Lichnerowicz Laplacian $\Delta_j:\Gamma(\Sym_0^j)\to \Gamma(\Sym_0^j)$ satisfy
\begin{gather}
\prod_{s=0}^j\left((T_j^+)^{\ast}T_j^+-a(s;j)\left(\Delta_j-b(s;j)c\right)\right)=0,
\label{eqn:fact2}\\ 
\textrm{where}\quad a(s;j)=\frac{(j-s+1)(n+j+s-2)}{(j+1)(n+2j-2)},\quad b(s;j)=j(n+j-1)+s(n+s-3).\nonumber
\end{gather}
\end{theorem}
\begin{proof}
We prove theorem by induction for $j$. The case of $j=0$ follows from the definition of $\Delta_0$, that is, $\Delta_0=(T_0^+)^{\ast}T_0^+$. We assume that the equation \eqref{eqn:fact2} holds for $j$. By \eqref{eqn:W-T}, 
\[
\begin{split}
 &T_j^+((T_j^+)^{\ast}T_j^+-a(s;j)(\Delta_j-b(s;j)c))\\
=&(T_j^+(T_j^+)^{\ast}-a(s;j)(\Delta_{j+1}-b(s;j)c))T_j^+\\
=&\underbrace{\frac{(j+2)(n+2j)}{(j+1)(n+2j-2)}}_{=:k(j)}((T_{j+1}^+)^{\ast}T_{j+1}^+-a(s;j+1)(\Delta_{j+1}-b(s;j+1)c))T_j^+.
\end{split}
\]
Then
\[
\begin{split}
0=&\prod_{s=0}^jT_j^+\left((T_j^+)^{\ast}T_j^+-a(s;j)\left(\Delta_j-b(s;j)c\right)\right)T_{j+1}^-\\
=&\prod_{s=0}^jk(j)\left((T_{j+1}^+)^{\ast}T_{j+1}^+-a(s;j+1)\left(\Delta_{j+1}-b(s;j+1)c\right)\right)T_j^+T_{j+1}^-.
\end{split}
\]
Finally we reach \eqref{eqn:fact2} for $j+1$. 
\end{proof}
If a tensor filed $\phi$ in $\Gamma(\Sym^j_0)$ satisfies $T_j^+\phi=0$, then we call it {\it trace-free conformal Killing tensor}, which has been investigated not only in differential geometry but also in relativity theory. As a corollary, we have an information of the eigenvalues of the Laplacian on trace-free conformal Killing tensors on a compact space of constant curvature (\cite{ST}, \cite{Take} for the sphere case). We introduce a filtration on $\Gamma(\Sym_0^j)$ as 
\[
F_j\subset F_{j+1}\subset \cdots \subset F_0=\Gamma(\Sym_0^j),\quad F_k=\ker T_{k}^-\cdots T_j^-.
\]
and an associated grading 
\[
\Gamma(\Sym_0^j)=\bigoplus_{0\le s\le j}W_s,\quad W_s:=T_{j-1}^+\cdots T_s^+(\ker T_s^-).
\]
\begin{corollary} On a compact space of constant curvature $K=c$, the eigenvalue of $\Delta_j$ on $\ker T_j^+\cap W_k$ is $(j(n+j-1)+k(n+k-3))c$. 
\end{corollary}

\section{Harmonic analysis for trace-free symmetric tensor fields on $S^n$}\label{sec:5}
The method to calculate the eigenvalues of the operators on the $n$-dimensional sphere $S^n=\SO(n+1)/\SO(n)$ is same to the spinor case. From Frobenius reciprocity and the branching rule, the space of trace-free symmetric tensor fields are decomposed with respect to $\SO(n+1)$, 
\[
L^2(S^n,\Sym_0^j)=\bigoplus_{0\le s\le j} \mathbf{V}_j(s)=\bigoplus_{0\le s\le j}\left(\bigoplus_{k\in \mathbb{Z}_{\ge 0}} V_j(k,s)\right).
\]
As mentioned in Section \ref{sec:3}, the space $V_j(k,s)$ stands for an irreducible $\SO(n+1)$-module with the highest weight $(k+j,s,0_{m-2})$ and its dimension is
\[
\frac{(2l+n-1)(2s+n-3)(l+n-2+s)(l+1-s)}{(n-1)(n-2)(s+n-3)(l+n-2)}\binom{l+n-2}{l+1}\binom{s+n-3}{s}
\]
by setting $l=k+j$. Note that, for $n=3$, $V_j(k,s)$ is the sum of $\SO(4)$-modules with the highest weight $(k+j,s)$ and $(k+j,-s)$. We also put $\mathbf{V}_j(s):=\bigoplus_{k\in \mathbb{Z}_{\ge 0}} V_j(k,s)$. As in \cite{Br1} Freudenthal's formula gives the eigenvalue of $\Delta_j$ on $V_j(k,s)$, 
\[
\left.\Delta_j\right|_{V_j(k,s)}=(j+k)(n+k+j-1)+s(s+n-3),
\]
so that we have 
\[
\left. \left(\Delta_j-b(s;j)\right)\right|_{V_j(k,s)}=k(n+k+2j-1).
\]
The factorization formula \eqref{eqn:fact2} allows us to guess the eigenvalue of $(T_j^+)^{\ast}T_j^+$ on $V_j(k,s)$ would be 
\begin{equation}
a(s;j)k(n+k+2j-1)=\frac{k(n+k+2j-1)(j-s+1)(n+j+s-2)}{(j+1)(n+2j-2)}.\label{eqn:eigen}
\end{equation}
We shall prove the claim that the eigenvalue of $(T_j^+)^{\ast}T_j^+$ on $V_j(k,s)$ coincides with \eqref{eqn:eigen} by induction for $j$. When $j=0$, the claim is true because of $\Delta_0=(T_0^+)^{\ast}T_0^+$. We suppose the claim holds for $j$. Then 
\[
\ker T_j^+=\bigoplus_{0\le s\le j}V_j(0,s)
\]
and $T_j^+:V_j(k,s)\to V_{j+1}(k-1,s)$ is isomorphism for any $k\ge 1$ and $s=0,\cdots,j$. For $\phi$ in $V_j(k,s)$, 
\[
\begin{split}
 &(j+2)(T_{j+1}^+)^{\ast}T_{j+1}^+(T_j^+\phi)\\
= &\left(\Delta_{j+1}+(n+j-2)(T_{j+1}^-)^{\ast}T_{j+1}^--2((j+1)(n+(j+1)-2))\right)(T_j^+\phi)\\
=&\left(\Delta_{j+1}+\frac{(j+1)(n+2j-2)}{n+2j}T_{j}^+(T_{j}^+)^{\ast}-2((j+1)(n+j-1))\right)(T_j^+\phi)\\
=&T_{j}^+\Delta_j\phi+\frac{k(n+k+2j-1)(j-s+1)(n+j+s-2)}{n+2j}T_j^+\phi\\  &\quad -2((j+1)(n+j-1))(T_j^+\phi)\\
=&\frac{(k-1)(n+k+2j)(j-s+2)(n+j+s-1)}{n+2j}T_j^+\phi.
\end{split}
\]
Thus, the eigenvalues of $(T_{j+1}^+)^{\ast}T_{j+1}^+$ on $V_{j+1}(k-1,s)$ for any $k\ge 1$ and $s=0,\cdots,j$ coincide with \eqref{eqn:eigen}. The cases of $V_{j+1}(k,j+1)$ for any $k\ge 0$ remain to be proved. Since $T_j^+$ is an overdetermined elliptic operator, the space of the sections $\Gamma(\Sym_0^{j+1})$ decomposes into the orthogonal direct sum of $\ker T_{j+1}^-$ and $\mathrm{Image}\; T_j^+$. From the above discussion, it holds that $\mathrm{Image}\; T_j^+$ is given by $\oplus_{0\le s\le j}\mathbf{V}_{j+1}(s)$. Then its orthogonal complement $\ker T_{j+1}^-$ is $\mathbf{V}_{j+1}(j+1)=\oplus_{k\ge 0}V_{j+1}(k,j+1)$. For $\phi$ in $V_{j+1}(k,j+1)$, by Weitzenb\"ock formula, we have
\[
\begin{split}
(T_{j+1}^+)^{\ast}T_{j+1}^+\phi=\frac{1}{j+2}(\Delta_{j+1}\phi-2((j+1)(n+(j+1)-2))\phi)
=\frac{k(n+k+2j+1)}{j+2}\phi,
\end{split}
\]
and check that this gives \eqref{eqn:eigen} on $V_{j+1}(k,j+1)$. Thus we have proved that \eqref{eqn:eigen} gives the eigenvalue of $(T_j^+)^{\ast}T_j^+$ on $V_j(k,s)$ for any $j$. We can also have the eigenvalue of the other operators by using Weitzenb\"ock formulas in Proposition \ref{prop:W-sym}. 
\begin{theorem}\label{thm:spec-sym}
On the $n$-dimensional standard sphere $S^n$, the eigenvalues of $(T_{j}^+)^{\ast}T_{j}^+$, $(T_{j}^-)^{\ast}T_{j}^-$, $U_j^{\ast}U_j$ and $\Delta_j$ on $V_j(k,s)$ are given by
\begin{gather}
\frac{k(n+k+2j-1)(j-s+1)(n+j+s-2)}{(j+1)(n+2j-2)},\\
\frac{(j-s)(k+1)(n+k+2j-2)(n+j+s-3)}{(n+j-3)(n+2j-2)},\\
\frac{s(k+j+1)(n+s-3)(n+k+j-2)}{(j+1)(n+j-3)},\\
\textrm{and}\quad (k+j)(k+j+n-1)+s(n+s-3),
\end{gather}
respectively. The kernels and images of the operators are realized as $\SO(n+1)$-modules in $L^2(S^n,\Sym_0^j)$ as follows:
\begin{align*}
\ker T_j^+=\bigoplus_{0\le s \le j}V_j(0,s),\quad \ker T_j^-=\mathbf{V}_j(j), \quad \ker U_j^-=\mathbf{V}_j(0),\\
\mathrm{Image}\;T_{j-1}^+=\bigoplus_{0\le s\le j-1}\mathbf{V}_j(s),\quad \mathrm{Image}\; T_{j+1}^-=L^2(S^n,\Sym_0^j)\ominus\ker T_j^+.
\end{align*}We also know that
\[
\ker T_j^+\cap \ker T_j^-=V_j(0,j),\quad \ker T_j^+\cap \ker U_j=V_j(0,0).
\]
\end{theorem}

As an interesting application to geometry, we shall discuss the space of Killing tensor fields on $S^n$ from the viewpoint of representation theory. We consider the differential $d$ on the symmetric tensor fields defined by 
\[
d:\Gamma(\Sym^j)\ni K\mapsto  dK:=\sum_{i=1}^ne_i\cdot \nabla_{e_i}K\in \Gamma(\Sym^{j+1}),
\]
where $\{e_i\}_i$ is a local orthonormal frame for $TM$ and $e_i\cdot$ denotes the symmetric tensor product by $e_i$. If $dK$ is zero, then $K$ is said to be a {\it Killing tensor field}. Since Killing tensor fields give the first integrals for geodesics, they play an important role in physics literature, especially in relativity theory. We refer to  \cite{HMS} for general results of Killing and conformal Killing tensor in Riemannian geometry. Since the differential $d$ is a derivation on $\oplus_j\Gamma(\Sym^j)$, that is, 
\[
d(K\cdot K')=(dK)\cdot K'+K\cdot (dK')\quad\textrm{ for $K,K'$ in $\oplus_j\Gamma(\Sym^j)$}, 
\]
the space $K(M)$ of the Killing tensor fields is a graded algebra, $K(M)=\oplus_j K^j(M)$, where $K^j(M)$ is the space of the Killing tensor fields with degree $j$. Let $K$ be a symmetric tensor fields with degree $j$. This $K$ decomposes  as
\[
K=K_0+g\cdot K_1+g^2\cdot K_2+\dots +g^l\cdot K_l \qquad (K_i\in \Gamma(\Sym_0^{j-2i}),\quad 0\le i\le l=[j/2])
\]
with respect to \eqref{symj}. Then it is easily shown that $K$ is a Killing tensor filed if and only if $\{K_i\}_{0\le i\le l}$ satisfy
\begin{equation}
(dK_0)_0 =0,\quad (dK_0)_1+g\cdot (dK_1)_0=0, \quad \dots, \quad (dK_{l-1})_1+g\cdot (dK_l)_0=0, \quad  (dK_l)_1=0. \label{eqn:killing}
\end{equation}
Note that if $K$ satisfies only the first condition $(dK_0)_0=0$, then $K$ is said to be a {\it conformal Killing tensor field}. For a trace-free symmetric tensor $K=K_0$ with degree $j$, $K$ is a trace-free Killing tensor (resp. trace-free conformal Killing tensor) if and only if  $K$ is in $\ker T_j^+\cap \ker T_j^-$ (resp. in $\ker T_j^+$). An important observation is that there is a nonzero constant $c=c(i,j)$ such that $(dK_i)_1+g\cdot (dK_{i+1})_0=0$ can be rewritten as $T_{j-2i}^-(K_i)=cT_{j-2i-2}^+(K_{i+1})$.

Now we consider the Killing tensor fields on the standard sphere. Let $K=K_0+g\cdot K_1+\cdots$ be a {\it primitive Killing tensor field} with degree $j$. In other words, $K$ is in $K^j(S^n)$, but is not the form of $K=K_0+g\cdot \tilde{K_1}$ with nonzero $\tilde{K}_1$ in $K^{j-2}(S^n)$. Then we may assume $K_0$ is in $V_j(0,s)$ because of $T_j^+(K_0)=0$. When $s=j$, we know $T_j^-(K_0)$ is zero by Theorem \ref{thm:spec-sym} and hence $T_{j-2}^+(K_1)$ is zero by \eqref{eqn:killing}. Thus we know $K$ is the form of $K=K_0+g\cdot \tilde{K_1}$ with $\tilde{K}_1$ in $K^{j-2}(S^n)$. Therefore $K=K_0$ due to primitiveness of $K$. When $s=j-1$, we show $T_j^-(K_0)=cT_{j-2}^+(K_1)$ is not zero and $K_1$ has to be in a $\SO(n+1)$-module with the highest weight $(j,j-1,0_{m-2})$ because $T_j^{\pm}$ is an invariant operator. Since $\Gamma(\Sym_0^{j-2})$ doesn't include such a component, there is no Killing tensor field $K$ with $K_0$ in $V_j(0,j-1)$. In the same manner, when $K_0$ is in $V_j(0,j-2i)$, the primitive Killing tensor $K$ is a form of $K=K_0+g\cdot K_1+\cdots+g^i\cdot K_i$ with $K_s$ in $V_{j-2s}(2s,j-2i)$ for $s=0,\cdots,i$. We also know there is no Killing tensor field $K$ with nonzero $K_0$ in $V_j(0,j-2i+1)$. 
\begin{proposition}[\cite{Take}]
Let $P^j(S^n)$ be the space of the primitive Killing tensors with degree $j$ on the standard sphere. Then 
\begin{gather*}
P^j(S^n)\cong \bigoplus_{0\le i\le [j/2]} V_j(0,j-2i)\cong \bigoplus_{0\le i\le [j/2]} (j,j-2i,0_{m-2}),\\
K^j(S^n)=\bigoplus_{0\le i\le [j/2]}g^i\cdot P^{j-2i}(S^n).
\end{gather*}

\end{proposition} 

\section{Higher spin Dirac operators on spinor fields with differential forms}\label{sec:6}

First we survey analysis for the generalized gradients on differential forms. Let $(M,g)$ be an oriented Riemannian manifold and $\Lambda^j(T^{\ast}M)$ be the bundle of differential forms. Due to the star operator $\ast:\Lambda^j(T^{\ast}M)\to \Lambda^{n-j}(T^{\ast}M)$, we have to study only the case of $j\le [n/2]$. As presented in \cite{H1}, \cite{SemC}, there are three differential operators 
\begin{gather*}
C:\Omega^j(M)=\Gamma(\Lambda^j(T^{\ast}M))\to \Gamma(\Lambda^{j,1}(T^{\ast}M)),\quad
 d:\Omega^j(M) \to \Omega^{j+1}(M),\quad d^{\ast}:\Omega^j(M) \to \Omega^{j-1}(M)
\end{gather*}
satisfying
\begin{align}
\Delta_j&=\nabla^{\ast}\nabla+\frac{1}{2}R_{\Lambda^j}=d^{\ast}d+dd^{\ast}, \quad  dd=0,\quad d^{\ast}d^{\ast}=0, \label{eqn:diff1}\\
\nabla^{\ast}\nabla&=C^{\ast}C+\frac{1}{j+1}d^{\ast}d+\frac{1}{n-j+1}dd^{\ast},\label{eqn:diff2}
\end{align}
where the highest weight of the bundle $\Lambda^{j,1}(T^{\ast}M)$ with respect to $\SO(n)$ is $(2,1_{j-1},0_{m-1-j})$ for $n=2m-1$ or $(2,1_{j-1},0_{m-j})$ for $n=2m$. Because of the equation \eqref{eqn:diff1} and ellipticity of $\Delta_j$, the Hodge-de Rham decomposition holds on a compact manifold $(M,g)$,
\[
\Omega^j(M)=H^j(M)\oplus d(\Omega^{j-1}(M))\oplus d^{\ast}(\Omega^{j+1}(M)).
\]
Therefore, if we know the eigenvalues of the Laplacian $\Delta_j$, then we have the eigenvalues of $d^{\ast}d$,  $dd^{\ast}$. In addition, from the action of curvature $R_{\Lambda^j}/2$, we can get the eigenvalues of $C^{\ast}C$. For example, on a space of constant curvature $K=c$, the curvature $R_{\Lambda^j}/2$ acts by a constant $j(n-j)c$. We will calculate the eigenvalues of the operators on the standard sphere in the next section. 

We move on to the higher spin Dirac operator on spinor fields coupled with differential forms. In other words, we shall study the case of the representation $\delta_j$ on $W_{\delta_j}$ with the highest weight 
\[
\begin{cases}
((3/2)_j,(1/2)_{m-j})\oplus((3/2)_j,(1/2)_{m-j-1},-1/2) & \textrm{for $n=2m$},\\
  ((3/2)_j,(1/2)_{m-1-j})   &\textrm{for $n=2m-1$},
   \end{cases}
\]
and the dimension $2^{[n/2]}\frac{n-2j+1}{n-j+1}\binom{n}{j}$, where $\delta_j$ splits into the sum of irreducible representation $\delta_j^{\pm}$ for $n=2m$. These spinor fields are realized as sections of the tensor bundle $S_0\otimes \Lambda^j(T^{\ast}M^c)$, that is, {\it spinor fields coupled with differential forms} (see Remark \ref{rem:spin-diff}). We consider the vector bundle with fiber $W_{\delta_j}$, 
\[
E_j:=\Spin(M)\times_{\Spin(n)}W_{\delta_j}.
\]
There are four generalized gradients on $\Gamma(E_j)$ for $0\le j\le [n/2]$,
\begin{align*}
& U_j:\Gamma(E_j)\to \Gamma(E_{j,1}) & \textrm{the (first) twistor operator},\\
 & \tilde{T}_j^+:\Gamma(E_j)\to \Gamma(E_{j+1}) & \textrm{the (second) twistor operator},\\
 & \tilde{D}_j:\Gamma(E_j)\to \Gamma(E_j) & \textrm{the higher spin Dirac opeator},\\
 & \tilde{T}_j^-:\Gamma(E_j)\to \Gamma(E_{j-1}) & \textrm{the co-twistor opeator},
\end{align*}
where $E_{j,1}$ is the vector bundle with the highest weight
\[
\begin{cases}
(5/2,(3/2)_{j-1},(1/2)_{m-j})\oplus (5/2,(3/2)_{j-1},(1/2)_{m-j-1},-1/2) & \textrm{for $n=2m$},\\
(5/2,(3/2)_{j-1},(1/2)_{m-1-j})& \textrm{for $n=2m-1$}.
\end{cases}
\]
\begin{rem}
When $j=0$, there are two generalized gradients on $E_0$, 
\[
\tilde{T}_0^+:\Gamma(E_0)\to \Gamma(E_{1}), \quad \tilde{D}_0:\Gamma(E_0)\to \Gamma(E_{0}).
\]
Then we set $U_0=0$ and $\tilde{T}_0^-=0$. When $n=2m$ and $j=m$, there are only two generalized gradients on the bundle $E_m=E_m^+\oplus E_m^-$ with the highest weight $((3/2)_{m})\oplus ((3/2)_{m-1},-3/2)$,
\[
U_m:\Gamma(E_m)\to \Gamma(E_{m,1}), \quad \tilde{T}_m^-:\Gamma(E_m)\to \Gamma(E_{m-1}).
\]
Note that the higher spin Dirac operator doesn't exist. Then we set $\tilde{D}_m=0$ and $\tilde{T}_m^+=0$ for $n=2m$. When $n=2m-1$ and $j=m-1$, there are three generalized gradients on the bundle $E_{m-1}$ with the highest weight $((3/2)_{m-1})$,
\[
U_{m-1}:\Gamma(E_{m-1})\to \Gamma(E_{m-1,1}), \quad \tilde{T}_{m-1}^-:\Gamma(E_{m-1})\to \Gamma(E_{m-2}),\quad \tilde{D}_{m-1}:\Gamma(E_{m-1})\to \Gamma(E_{m-1}).
\]
Then we set $\tilde{T}_{m-1}^+=0$ for $n=2m-1$. 
\end{rem}

From now on, we assume that $(M,g)$ is a space of constant curvature $K=c$ with a spin structure. We calculate the Weitzenb\"ock formulas on $\Gamma(E_j)$ from \cite{H1}, 
\begin{equation}
\begin{split}
\Delta_j=&\nabla^{\ast}\nabla+\frac{1}{2}R_{\delta_j}=\nabla^{\ast}\nabla+(j(n-j+1)+\frac{n(n-1)}{8})c  \\
=&\frac{(n+2)(n-2j)}{n-2j+2}\tilde{D}_j^2 +\frac{4(n-2j+1)(n-j+2)}{(n-2j+3)(n-2j+2)}(\tilde{T}_{j}^-)^{\ast}\tilde{T}_j^-+ (j(n-j)-\frac{n(n-1)}{8})c \\
  =&-\frac{4(n-2j+1)(j+1)}{(n-2j-1)(n-2j)} (\tilde{T}_{j}^+)^{\ast}\tilde{T}_j^+        +\frac{(n-2j+2)(n+2)}{n-2j}\tilde{D}_j^2+((j-1)(n-j+1)-\frac{n(n-1)}{8})c \\
=& \frac{(n-j+2)(n+2)}{(n-j+1)(n+1)}U_j^{\ast}U_j+\frac{(n-2j)(n-2j+1)}{(n-2j-1)(n-j+1)}(\tilde{T}_{j}^+)^{\ast}\tilde{T}_j^+ +(j(n-j+2)+\frac{n(n+1)}{8})c.
\end{split} \label{eqn:estimate}
\end{equation}
If the denominator of the above equation has zero, then we ignore the equation. For example, when $n=2m$ and $j=m$, we don't consider the third equality. From some formulas in  \cite{H2} and the equation \eqref{eqn:rel} in Lemma \ref{relative}, we also have 
\begin{gather*}
\frac{1}{\sqrt{n-2j}}\tilde{D}_j\tilde{T}_{j-1}^+=\frac{1}{\sqrt{n-2j+4}}\tilde{T}_{j-1}^+\tilde{D}_{j-1},\quad \tilde{T}_{j+1}^+\tilde{T}_{j}^+=0,\quad \tilde{T}_{j-1}^-\tilde{T}_{j}^-=0,\\
(\tilde{T}_{j-1}^+)^{\ast}\tilde{T}_{j-1}^+=\frac{(n-2j+1)(n-j+2)}{j(n-2j+3)}\tilde{T}_j^-(\tilde{T}_j^-)^{\ast}. 
\end{gather*}
Then we normalize the operators as follows:
\begin{gather*}
D_j:=\sqrt{\frac{(n+2)(n-2j)}{n-2j+2}}\tilde{D}_j,\\
T_j^-:=2\sqrt{\frac{(n-2j+1)(n-j+2)}{(n-2j+3)(n-2j+2)}}\tilde{T}_j^-,\quad T_j^+:=2\sqrt{\frac{j+1}{n-2j}}\tilde{T}_j^+. 
\end{gather*}
\begin{theorem}\label{thm:bwsd}
The generalized gradients on $\Gamma(E_j)$ satisfy 
\begin{gather*}
\begin{split}
\Delta_j=&D_j^2+(T_{j}^-)^{\ast}T_j^-+ (j(n-j)-\frac{n(n-1)}{8})c\\
  =&-\frac{n-2j+1}{n-2j-1} (T_{j}^+)^{\ast}T_j^+ +\frac{(n-2j+2)^2}{(n-2j)^2}D_j^2 +((j-1)(n-j+1)-\frac{n(n-1)}{8})c,
\end{split}\\
D_jT_{j-1}^+=\frac{n-2j}{n-2j+2}T_{j-1}^+D_{j-1},\quad T_{j+1}^+T_j^+=0,\quad T_{j-1}^-T_j^-=0,\quad T_{j-1}^+=(T_j^-)^{\ast}.
\end{gather*}
\end{theorem}

\begin{rem}\label{rem:spin-diff}
In \cite{BS}, these operators are realized  as components of the twisted Dirac operator $D(j)$ on $S_0\otimes \Lambda^j(T^{\ast}M^c)$ up to a constant multiple. This operator $D(j)$ decomposes along $S_0\otimes \Lambda^j(T^{\ast}M^c)=\oplus_{k=0}^j E_k$, 
\[
D(j)=
\begin{pmatrix}
D(j)_j & T(j)_{j-1}^+   & 0       & 0& \dots & 0\\
T(j)_j^- & D(j)_{j-1} & T(j)_{j-2}^+ & 0 & \cdots & 0\\
0 &  T(j)_{j-1}^- & D(j)_{j-2}' &  T(j)_{j-3}& \ddots & 0\\
\vdots  & \ddots   & \ddots & \ddots & \ddots   & \vdots\\
0   &  0&  0& \dots           & D(j)_1        & T(j)^+_0\\
 0   &  0&  0& \dots           & T(j)_1^-        & D(j)_0
\end{pmatrix}
\]
where $D(j)_k$ (resp. $T(j)_k^{\pm}$) is a nonzero constant multiple of  $\tilde{D}_k:\Gamma(E_k)\to \Gamma(E_k)$ (resp. $\tilde{T}_k^{\pm}:\Gamma(E_k)\to \Gamma(E_{k\pm 1})$) for $k=0,\dots, j$. Considering the square of $D(j)$, we have some formulas equivalent to Theorem \ref{thm:bwsd} and check that $D(j)_j$, $T(j)_j^-$ and $T(j)_{j-1}^+$ coincide with our $D_j$, $T_j^-$ and $T_{j-1}^+$, respectively. 
\end{rem}
As applications of theorem, we obtain a factorization formula and the Hodge-de Rham decomposition with respect to the operators on $\Gamma(E_j)$. 
\begin{corollary}[factorization formula]
We consider a space $(M,g)$ of constant curvature $K=c$ with a spin structure. On $\Gamma(E_j)$ for $1\le j\le[n/2]-1$, we have
\[
\begin{split}
0=\left(D_j^2-\frac{(n-2j)^2}{(n-2j+2)^2} (\Delta_j-((j-1)(n-j+1)-\frac{n(n-1)}{8})c )\right)\\ \times \left(D_j^2-(\Delta_j- (j(n-j)-\frac{n(n-1)}{8})c)\right).
\end{split}
\]
\end{corollary}
\begin{proof}
By $(T_j^+)^{\ast}T_j^+ (T_j^-)^{\ast}T_j^-=(T_j^+)^{\ast}T_j^+ T_{j-1}^+T_j^-=0$, we can prove the corollary. 
\end{proof}
Deleting $D_j^2$ from two equations for $\Delta_j$ in Theorem \ref{thm:bwsd}, we get 
\begin{align*}
\Delta_j=\frac{(n-2j)^2}{4(n-2j-1)}(T_{j}^+)^{\ast}T_j^+ +\frac{(n-2j+2)^2}{4(n-2j+1)}(T_{j}^-)^{\ast}T_j^- +\frac{n(n+1)}{8}c.
\end{align*}
Taking the usual Hodge-de Rham decomposition as a model, we can show the following $L^2$-orthogonal decomposition.  
\begin{proposition}
On a compact spin manifold $(M,g)$ of constant sectional curvature $K=c$, we have the Hodge-de Rham decomposition for spinor fields coupled with differential forms,
\[
\begin{split}
\Gamma(E_j)=&T_{j-1}^+(\Gamma(E_{j-1}))\oplus T_{j+1}^-(\Gamma(E_{j+1}))\oplus \ker (\Delta_j-\frac{n(n+1)}{8}c),\\
\ker  T_j^+=&T_{j-1}^+(\Gamma(E_{j-1}))\oplus \ker (\Delta_j-\frac{n(n+1)}{8}c),\\
\ker  T_j^-=&T_{j+1}^-(\Gamma(E_{j-1}))\oplus \ker (\Delta_j-\frac{n(n+1)}{8}c).\end{split}
\]
Here we exclude the case of $n=2m-1$ and $j=m-1$. When $j=0$, the first equation means
\[
\Gamma(E_0)=T_{1}^-(\Gamma(E_{1}))\oplus \ker (\Delta_{0}-\frac{n(n+1)}{8}c).
\]
When $n=2m$ and $j=m$, it does
\[
\Gamma(E_m)=T_{m-1}^+(\Gamma(E_{m-1}))\oplus \ker (\Delta_m-\frac{n(n+1)}{8}c).\]
\end{proposition}
Note that $\ker (\Delta_j-\frac{n(n+1)}{8}c)$ is zero on a compact spin manifold of positive constant curvature for $j\ge 1$ because we have $\Delta_j-\frac{n(n+1)}{8}c\ge  j(n-j+2)c$ from the last equality in \eqref{eqn:estimate}.

\section{Harmonic analysis for spinor fields with differential forms on $S^n$}\label{sec:7}
We review a well-known result for the eigenvalues of the operators acting on differential forms on the standard sphere, \cite{Folland}, \cite{IT}. Let $V_j(k)$ be an irreducible $\SO(n+1)$-module with the highest weight
\[
(k+1,1_{j-1},0_{m-j})=(k,0_{m-1})+(1_j,0_{m-j})
\]
for $0\le j\le m=[\frac{n+1}{2}]$. Then Frobenius reciprocity and branching rule give us the decomposition of $L^2(S^n,\Lambda^j(T^{\ast}M^c))$ as an $\SO(n+1)$-module.  
\begin{enumerate}
\item For $j=0$,
\[
L^2(S^{n},\Lambda^0(T^{\ast}M^c))\cong  \bigoplus_{k\ge 0}V_{1}(k)\oplus V_0(0),\quad \ker d=H^0(S^n)\cong V_{0}(0).
\]
\item For $n=2m$ and $j=m$,
\[
L^2(S^{2m},\Lambda^m(T^{\ast}M^c))\cong \bigoplus_{k\ge 0}2V_{m}(k),\quad 
\ker d^{\ast}\cong \bigoplus_{k\ge 0}V_{m}(k),\quad \ker d\cong \bigoplus_{k\ge 0}V_{m}(k).
\]
\item For $n=2m-1$ and $j=m-1$,
\[
\begin{split}
L^2(S^{2m-1},\Lambda^{m-1}(T^{\ast}M^c))&\cong  \bigoplus_{k\ge 0} V_{m}(k)\oplus \bigoplus_{k\ge 0}V_{m-1}(k),\quad   V_m(k):=V_{m}^+(k)\oplus V_{m}^-(k),\\
\ker d^{\ast}&\cong \bigoplus_{k\ge 0}V_{m}(k),\quad \ker d\cong \bigoplus_{k\ge 0}V_{m-1}(k),
\end{split}
\]
where $V_m^{\pm}(k)$ is the irreducible $\SO(2m)$-module with the highest weight $(k+1,1_{m-2},\pm 1)$. 
\item Otherwise,
\[
L^2(S^{n},\Lambda^j(T^{\ast}M^c))\cong \bigoplus_{k\ge 0}V_{j+1}(k)\oplus \bigoplus_{k\ge 0}V_{j}(k),\quad 
\ker d^{\ast}\cong \bigoplus_{k\ge 0}V_{j+1}(k),\quad \ker d\cong \bigoplus_{k\ge 0}V_{j}(k).
\]
\end{enumerate}
We calculate the eigenvalue of $\Delta_j$ on $V_j(k)$ by Freudenthal's formula. Then the next proposition follows from \eqref{eqn:diff1} and \eqref{eqn:diff2}.
\begin{proposition}[\cite{Folland}, \cite{IT}]
The eigenvalues of the operators on $V_{j+1}(k)$ and $V_j(k)$ are given as follows:
\[
\begin{array}{|c|c|c|}
\hline 
      & \textrm{on}\; V_{j+1}(k)\subset \ker d^{\ast}  & \textrm{on}\; V_j(k) \subset \ker d\\
\hline
\Delta_j &  (k+j+1)(n+k-j)  &  (k+j)(n-j+k+1)  \\
\hline
 dd^{\ast} & 0 & (k+j)(n-j+k+1)\\
\hline
d^{\ast}d & (k+j+1)(n+k-j) & 0\\
\hline
 C^{\ast}C & \frac{j}{j+1}k(n+k+1) & \frac{n-j}{n-j+1}k(n+k+1) \\
\hline
\end{array}
\]
In particular, the space of the Killing $j$-forms $\ker C\cap \ker d^{\ast}$ is isomorphic to $V_{j+1}(0)$, and the space of the co-Killing $j$-forms $\ker C\cap \ker d$ is isomorphic to $V_{j}(0)$. 
\end{proposition}
The (co-) Killing forms have interesting geometric meaning like Killing tensor fields \cite{SemC}.  

We shall study harmonic analysis for the spinor fields coupled with differential forms on the standard sphere. The space $L^2(S^n,E_j)$ is decomposed as $\Spin(n+1)$-module. 
\begin{proposition}
\begin{enumerate}
\item For the case of $n=2m$, we denote by $V_j(k)'$ an irreducible $\Spin(2m+1)$-module with the highest weight
\[
(k+3/2,(3/2)_{j-1},(1/2)_{m-j})=(k+1/2,(1/2)_{m-1})+(1_j,0_{m-j})
\]
for $k=0,1,2,\dots$. Then
\begin{enumerate}
\item For $j=0$, 
\begin{gather*}
L^2(S^{2m},E_0^{\pm})\cong\bigoplus_{k\ge 0} V_{1}(k)'\oplus V_0(0)', \quad  L^2(S^{2m},E_0)\cong\bigoplus_{k\ge 0} 2V_{1}(k)'\oplus  2V_0(0)', \\
\ker T_0^+=2V_0(0)'.
\end{gather*}
\item For $1\le j\le m-1$,
\begin{gather*}
L^2(S^{2m},E_j^{\pm})\cong\bigoplus_{k\ge 0} V_{j+1}(k)'\oplus \bigoplus_{k\ge 0}V_j(k)', \quad  L^2(S^{2m},E_j)\cong\bigoplus_{k\ge 0} 2V_{j+1}(k)'\oplus \bigoplus_{k\ge 0}2V_j(k)', \\
\ker T_j^-=\bigoplus_{k\ge 0} 2V_{j+1}(k)', \quad \ker T_j^+=\bigoplus_{k\ge 0} 2V_{j}(k)'.
\end{gather*}
\item For $j=m$,
\begin{gather*}
L^2(S^{2m},E_m^{\pm})\cong\bigoplus_{k\ge 0} V_{m}(k)', \quad  L^2(S^{2m},E_j)\cong\bigoplus_{k\ge 0} 2V_{m}(k)', \quad \ker T_m^-=\{0\}.
\end{gather*}
\end{enumerate}
\item For the case of $n=2m-1$, we denote by $V_j^{\pm}(k)'$ an irreducible $\Spin(2m)$-module with the highest weight
\[
(k+3/2,(3/2)_{j-1},(1/2)_{m-j-1},\pm 1/2)=(k+1/2,(1/2)_{m-1},\pm 1/2)+(1_j,0_{m-j})
\]
for $k=0,1,2,\dots$, and put $V_j(k)'=V_j^+(k)\oplus V_j^-(k)$. Then
\begin{enumerate}
\item For $j=0$,
\begin{gather*}
L^2(S^{2m-1},E_0)\cong\bigoplus_{k\ge 0} V_{1}(k)'\oplus V_0(0)', \quad \ker T_0^+=V_0(0)'.
\end{gather*}
\item For $1\le j\le m-2$,
\begin{align*}
L^2(S^{2m-1},E_j)\cong \bigoplus_{k\ge 0} V_{j+1}(k)'\oplus \bigoplus_{k\ge 0} V_{j}(k)',\quad 
\ker T_j^-=\bigoplus_{k\ge 0} V_{j+1}(k)', \quad \ker T_j^+=\bigoplus_{k\ge 0} V_{j}(k)'.
\end{align*}
\item For $j=m-1$,
\begin{align*}
L^2(S^{2m-1},E_{m-1})\cong \bigoplus_{k\ge 0} V_{m}(k)'\oplus \bigoplus_{k\ge 0} V_{m-1}(k)',\quad 
\ker T_{m-1}^-=\bigoplus_{k\ge 0} V_{m}(k)'.
\end{align*}

\end{enumerate}

\end{enumerate}
\end{proposition}
We have already known how the sections of $\{E_j\}_{j}$ relate to each other through operators $\{T_j^{\pm}\}_j$. In fact, $\ker T_j^+=\mathrm{Image}\; T_{j-1}^+$ and $\ker T_j^-=\mathrm{Image} \; T_{j+1}^-$. Then all we have to do is to calculate the eigenvalues of the operators on $\Gamma(E_j)$. The eigenvalues of $D_j^2$ were calculated in \cite{BS} by using a method from a viewpoint of parabolic geometry in \cite{BOO}. We calculate them only from the eigenvalues of $\Delta_j$ and Weitzenb\"ock formulas. 
\begin{theorem}
The eigenvalues of $D_j^2$, $(T_j^-)^{\ast}T_j^-$, $(T_j^+)^{\ast}T_j^+$ and $U_j^{\ast}U_j$ on $L^2(S^n,E_j)$ are given as follows: 
\[
\begin{array}{|c|c|c|}
\hline 
      & \text{on}\; V_{j+1}(k)'\subset \ker T_j^-  & \text{on}\; V_j(k)' \subset \ker T_j^+\\
\hline
\Delta_j &  (k+j+1)(n-j+k+1)+\frac{n(n+1)}{8}    & (k+j)(n-j+k+2)+\frac{n(n+1)}{8} \\
\hline
D_j^2 &   \left(\frac{n}{2}+k+1\right)^2  &  \frac{(n-2j)^2}{(n-2j+2)^2}\left(\frac{n}{2}+k+1\right)^2 \\
\hline
 (T_j^-)^{\ast}T_j^- & 0 &  \frac{4(n-2j+1)}{(n-2j+2)^2}(k+j)(n-j+k+2)\\
\hline
 (T_j^+)^{\ast}T_j^+ & \frac{4(n-2j-1)}{(n-2j)^2}(n-j+k+1)(k+j+1)  &  0 \\
\hline
 U_j^{\ast}U_j & \frac{(n+1)j}{(n+2)(j+1)}k(n+k+2) & \frac{(n+1)(n-j+1)}{(n+2)(n-j+2)}k(n+k+2)\\
\hline
\end{array} 
\]
In particular, for $n=2m-1$ (resp. $n=2m$), $\ker U_j\cap \ker T_j^-=V_{j+1}(0)'$ (resp. $2V_{j+1}(0)'$) and $\ker U_j\cap \ker T_j^+=V_j(0)'$ (resp. $2V_{j}(0)'$). 
\end{theorem}

\appendix

\section{Calculation}\label{app:A}
In this appendix, we give calculation of the eigenvalues of Casimir elements and Weitzenb\"ock formulas needed in Section \ref{sec:2}. For more detail, see \cite{H1}. Let $(\pi_j,W_j)$ be an irreducible representation with the highest weight $(j+1/2,(1/2)_{[n/2]-1})$. The conformal weights $\{w(\lambda)\}_{\lambda}$ and the shifted conformal weights $\{\hat{w}(\lambda)=w(\lambda)+\frac{n-1}{2}\}_{\lambda}$ associated to the irreducible summand of $W_j\otimes \mathbb{C}^n$ are given by 
\begin{align*}
w(\pi_{j+1})&:=w(\pi_{j+1};\pi_j)=j+1/2,\quad \hat{w}(\pi_{j+1})=j+n/2,\\
w(\pi_{j,1})&:=w(\pi_{j,1};\pi_j)=-1/2,\quad \hat{w}(\pi_{j,1})=n/2-1,\\
w(\pi_{j})&:=w(\pi_{j};\pi_j)=-n/2+1/2,\quad \hat{w}(\pi_{j})=0,\\
w(\pi_{j-1})&:=w(\pi_{j-1};\pi_j)=-n-j+3/2,\quad \hat{w}(\pi_{j-1})=1-j-n/2.
\end{align*}
Then we can easily calculate the relative dimensions,
\begin{gather*}
\frac{\dim W_{j+1}}{\dim W_j}=\frac{n+j-1}{j+1},\quad \frac{\dim W_{j,1}}{\dim W_j}=\frac{(n-3)(n+j-1)j}{(n+j-2)(j+1)},\quad \frac{\dim W_{j}}{\dim W_j}=1,\quad
\frac{\dim W_{j-1}}{\dim W_j}=\frac{j}{n+j-2}.
\end{gather*}
The Casimir element $c_k$ with order $k$ is an $\SO(n)$-invariant element in the enveloping algebra $U(\mathfrak{so}(n))$ given by 
\[
c_k=\sum_{1\le i_1,i_2\dots i_{k}\le n} e_{i_1i_2}e_{i_2i_3}\dots e_{i_{k}i_1},
\]
where $\{e_{ij}=e_i\wedge e_j\}_{i<j}$ is a standard basis for $\mathfrak{so}(n)$. Shifting $e_{ij}$ to $\hat{e}_{ij}=e_{ij}+\frac{n-1}{2}$, we define the shifted Casimir element $\hat{c}_k=\sum \hat{e}_{i_1i_2}\hat{e}_{i_2i_3}\dots \hat{e}_{i_{k}i_1}$. The eigenvalue of $c_k$ on $W_j$ is computed by the conformal weights and the relative dimensions.
\begin{align*}
\pi_j(c_k)=\sum_{\lambda} w(\lambda)^k\frac{\dim W_{\lambda}}{\dim W_j},\quad 
\pi_j(\hat{c}_k)=\sum_{\lambda} \hat{w}(\lambda)^k\frac{\dim W_{\lambda}}{\dim W_j},
\end{align*}
where $\lambda$ runs among $\pi_{j+1}$, $\pi_{j,1}$, $\pi_{j}$ and $\pi_{j-1}$.  Since the number of the generalized gradients on $S_j$ is four, there are two independent Weitzenb\"ock formulas. In our case we have the following two independent formulas on $S_j$. One is 
\begin{equation*}
 -\frac{1}{2}R^1_j=w(\pi_{j+1})(\tilde{T}_j^+)^{\ast}\tilde{T}_j^++w(\pi_{j,1}) U_j^{\ast}U_j+ w(\pi_{j}) \tilde{D}_j^2+w(\pi_{j-1})(\tilde{T}_j^-)^{\ast}\tilde{T}_j^-.
\end{equation*}
The second one is complicated,
\begin{equation*}
\hat{R}^4_j=a(\pi_{j+1})(\tilde{T}_j^+)^{\ast}\tilde{T}_j^+ +a(\pi_{j,1}) U_j^{\ast}U_j+ a(\pi_j)\tilde{D}_j^2+a(\pi_{j-1})(\tilde{T}_j^-)^{\ast}\tilde{T}_j^-,
\end{equation*}
where
\[
a(\lambda)=\sum_{p=0}^3\pi_j(\hat{c}_{3-p})(-\hat{w}(\lambda))^p,\quad \textrm{for $\lambda=\pi_{j+1},\pi_{j,1},\pi_{j},\pi_{j-1}$}.
\]
Furthermore, the curvature actions $R^1_j$ and $\hat{R}^4_j$ on the standard sphere are 
\[
R^1_j=\pi_j(c_2),\quad \hat{R}^4_j=\pi_j(\hat{c}_5-\frac{n-1}{2}\hat{c}_4).
\]
On the other hand, it follows form the definition of the generalized gradients that 
\[
\nabla^{\ast}\nabla=(\tilde{T}_j^+)^{\ast}\tilde{T}_j^++U_j^{\ast}U_j+\tilde{D}_j^2+(\tilde{T}_j^-)^{\ast}\tilde{T}_j^-.
\]
Deleting two operators from the above equation by using two Weitzenb\"ock formulas, we obtain \eqref{WF1} and \eqref{WF2}.


\end{document}